\documentclass[12pt]{article}

\usepackage{tipa}
\usepackage{amsfonts,amssymb}
\usepackage{amssymb}
\usepackage{graphicx}
\usepackage{amsthm}
\usepackage{amsmath}
\usepackage[hang]{caption2}
\usepackage[all]{xy}\xyoption{all}

\newtheorem{theorem}{Theorem}
\newtheorem{proposition}{Proposition}

\newtheorem{remark}{Remark}
\newtheorem{definition}{Definition}

\begin{document}

\title{Symmetric orbits arising from Figure-Eight for N-body
problem}

\author{Leshun Xu$^1$ and Yong Li$^2$\\
$^1$College of Mathematics and Computer Science,\\
Nanjing Normal University, Nanjing, 210097, P.R. China.\\
$^2$College of Mathematics,\\
Jilin University, Changchun, 130012, P.R. China.\\
e-mail: lsxu@sohu.com {\rm and} liyong@jlu.edu.cn }

\date{11.2006}

\maketitle

\begin{abstract}
In this paper, we first describe how we can arrange any bodies on
Figure-Eight without collision in a dense subset of $[0,T]$ after
showing that the self-intersections of Figure-Eight will not happen
in this subset. Then it is reasonable for us to consider the
existence of generalized solutions and non-collision solutions with
Mixed-symmetries or with Double-Eight constraints, arising from
Figure-Eight, for N-body problem. All of the orbits we found
numerically in Section \ref{se7} have not been obtained by other
authors as far as we know. To prove the existence of these new
periodic solutions, the variational approach and critical point
theory are applied to the classical N-body equations. And along the
line used in this paper, one can construct other symmetric
constraints on N-body problems and prove the existence of periodic
solutions for them.
\end{abstract}

MSC-class: 70F10, 37C80, 70G75

\section{Introduction}
N-body problem is a classical problem in mathematics and celestial
mechanics. Since Euler's collinear solutions and Lagrange's
equilateral triangle solutions, any more ``clearly'' (including
the existence and the shape of the orbit) solutions have not been
found except Figure-Eight solution, a planer three-body orbits
with equal masses, which was proved strictly by Chenciner and
Montgomery in \cite{Chen1} in 1999. This remarkable and
interesting orbit arouse many authors' curiosity in studying more
choreography orbits or other interesting symmetric constraints
acting on the orbits of N-body problems, see \cite{{SU6},{Ck},
{Chen2},{Su1},{RM1},{Sh},{S1},{S2}}, etc. Furthermore, some
investigations of a possible star system in real space have been
performed in practice (see \cite{no}).

Among those authors' studies, the variational approach and the
theory of critical points perform an important role. By minimizing
the Lagrangian action on a space of loops symmetry with respect to a
well-chosen symmetry group, which also has some simple conditions on
it, Ferrario and Terracini give a fairly general condition on
symmetry groups of the loop space in \cite{Su1}. But their results
``are not suitable to prove the existence of those orbits that are
found as minimizers in classes of paths characterized by homotopy
conditions or a mixture of symmetry and homotopy conditions''. A
further study based on \cite{Su1} is carried by Shibayama to
distinguish the solutions, which have been obtained by Shibayama in
\cite{Sh} and some of which by Chen in \cite{Ck}.

In this paper, we study the Figure-Eight constraint carefully in
Section \ref{se3}, and draw a conclusion that even bodies with
equivariant shift on Figure-Eight will collide at least on the
original point, and odd bodies will not collide in a dense subset of
$[0,T]$. Then in Section \ref{se4}, we compose Figure-Eight from
simple to double (or more, that we do not describe in this paper) by
setting odd bodies on them reasonably. And in Section \ref{se5}, it
is practical to mix Figure-Eight with other constraints on $q^{(3)}$
in space. Then along Gordon \cite{{Go1},{Go2}} and Zelati's
\cite{Z4} line, in Section \ref{se6}, we have proved the existence
of generalized solutions and non-linear solutions for these periodic
orbits. The numerical study in Section \ref{se4} highly enlightens
us that we can arrange many (even or odd) bodies on Figure-Eight by
choosing an appropriate {\em phase differences}. Obviously, all the
study is based on Figure-Eight constraint, hence this paper is a
continued study on Chenciner and Montgomery's significant work.

\section{N-body problem and Figure-Eight}
The periodic solutions for N-body problem is described by the
following second order ordinary differential equations in general:
\begin{equation}
\label{N} -m_{i}\ddot{q_{i}} =\nabla U(q),\quad i=1,\cdots,N,
\end{equation}
\begin{equation}
\label{P} q(0)=q(T), \quad \dot{q}(0)=\dot{q}(T),
\end{equation}
where $m_{i}$ is the mass of $i$th body, $q_{i}$ is the position of
$i$th body in $\mathbb{R}^{d}$, $N$ is the number of bodies, $T$ is
the period, and generally the singular potential function $U$ is
\begin{equation*}
U(q)=\frac{1}{2}\sum_{1\leq i\neq j\leq N}U_{ij}(q_{j}-q_{i}),
\end{equation*}
with assumptions: $\forall 1\leq i\neq j \leq N$ and $\forall \xi\in
R^{d}\backslash \{0\}$,
\begin{equation*}\label{a1}
U_{ij}\in C^{2}(R^{d}\backslash \{0\},R),\quad
U_{ij}(\xi)=U_{ji}(\xi),
\end{equation*}
\begin{equation*}\label{a2}
U_{ij}(\xi) \rightarrow -\infty,\quad\quad\quad\quad \mbox{as} \quad
|\xi|\rightarrow 0,
\end{equation*}

\begin{equation*}\label{a3}
U(q_{1},\cdots,q_{N})\leq 0, \quad\quad
\forall(q_{1},\cdots,q_{N})\in (R^{d})^{N},
\end{equation*}
where $U_{ij}$ is usually described by
\begin{equation}\label{Uil1}
U_{ij}=\frac{m_{i}m_{j}}{|q_{i}-q_{j}|^{\alpha}}.
\end{equation}
When $\alpha=1$ in (\ref{Uil1}), it is called {\em Newton} N-body
problem. The figures, as the numerical examples to illuminate the
existence of the solutions in this paper, are all coming from Newton
N-body problem. Poincar\'e has shown that it is impossible for
$N\geq 3$ to find an explicit expression for the general solution.
So even for the three-body problem, it is impossible to describe all
the solutions. Then the numerical approach and/or endowing orbits
with some symmetric constraints are quite reasonable directions for
the study in N-body problem. No matter which direction adopted by
the authors it is, the Lagrangian
\begin{equation*}
J:\mathbb{H}^{1}(\mathbb{R}/T\mathbb{Z},(\mathbb{R}^{d})^{N})\rightarrow\mathbb{R},
\end{equation*}
\begin{equation}\label{J}
J(q)=\frac{1}{2}\sum_{i=1}^{N}m_{i}\int_{0}^{T}|\dot{q}_{i}|^{2}dt-\int_{0}^{T}U(q)dt,
\end{equation}
acts an important role in their study in finding the periodic
solutions for N-body problems.

In one direction, finding some symmetric orbits, Chenciner and
Montgomery have introduced a new ``torsional'' symmetric constraint
to the possible orbits for Three-body problem with equal masses in
1999. The ``torsional'' symmetry is described by the following
expressions:

Let $q: (\mathbb{R}/T\mathbb{Z})\rightarrow \mathbb{R}^2$, $i=1, 2,
3$, the orbit of the $i$-th body is
\begin{equation}
\label{e1} q_{i}=q\left(t+(3-i)\frac{T}{3}\right),
\end{equation}

\begin{equation}
\label{e2} q(t)+q(t+\frac{T}{3})+q(t+\frac{2T}{3})=0, \quad
t\in[0,T],
\end{equation}

\begin{equation}
\label{e3} q(\sigma(t))=\sigma(q(t)),\quad q(\tau(t))=\tau(q(t)).
\end{equation}
where $q(t)=(q^{(1)}(t),q^{(2)}(t))$, $\sigma$ and $\tau$ are
generators of the action of Klein group
$\mathbb{Z}/2\mathbb{Z}\times\mathbb{Z}/2\mathbb{Z}$ on
$\mathbb{R}/T\mathbb{Z}$ and on $\mathbb{R}^{2}$:

\begin{equation*}
\sigma(t)=t+\frac{T}{2},\quad \tau(t)=-t+\frac{T}{2},
\end{equation*}
\begin{equation*}
\sigma(q)=(-q^{(1)},q^{(2)}),\quad \tau(q)=(q^{(1)},-q^{(2)}).
\end{equation*}

By the constraints (\ref{e1})-(\ref{e3}), the existence is proved in
\cite{Chen1} for three equal masses with the exclusion of collision.
It is called simple choreography orbits, Figure-Eight solution.
After that, many choreography/choreographies solutions are found
numerically by many authors, e.g. Sim\'{o} \cite{{S1},{S2}}, Arioli,
Barutello and Terracini \cite{SU6} etc.

\section{Number of the bodies on Figure-Eight}\label{se3}
One of the facile development to simple choreography is considering
more bodies on the same Figure-Eight. We can develop (\ref{e1}) and
(\ref{e2}) by
\begin{equation}\label{e4}
q_{i}=q\left(t+(N-i)\frac{T}{N}\right), \quad i=1,\cdots,N,
\end{equation}
\begin{equation}\label{e5}
\sum^{N}_{i=1}q_{i}(t)=0, \quad t\in[0,T].
\end{equation}
Hence, (\ref{e4}), (\ref{e5}) and (\ref{e3}) describe $N$ bodies
on the same Figure-Eight with equivalent shift $\frac{T}{N}$. By
(\ref{e3}), one can deduce that
\begin{equation}\label{q1}
\left(q^{(1)}(t),q^{(2)}(t)\right)=\left(-q^{(1)}(t+\frac{T}{2}),q^{(2)}(t+\frac{T}{2})\right),
\end{equation}
\begin{equation}\label{q2}
\left(q^{(1)}(t),q^{(2)}(t)\right)=\left(q^{(1)}(-t+\frac{T}{2}),-q^{(2)}(-t+\frac{T}{2})\right).
\end{equation}
When we set $t=-t$,
\begin{eqnarray*}
\left(q^{(1)}(-t),q^{(2)}(-t)\right)&=&\left(q^{(1)}(t+\frac{T}{2}),-q^{(2)}(t+\frac{T}{2})\right)\\
&=&-\left(q^{(1)}(t),q^{(2)}(t)\right),
\end{eqnarray*}
It shows that $q(t)$ is an odd function. Let $t=0$, one can deduce
the following formulas from (\ref{q1}) and (\ref{q2}),
\begin{equation*}
\left(-q^{(1)}(0),q^{(2)}(0)\right)=\left(q^{(1)}\left(\frac{T}{2}\right),q^{(2)}\left(\frac{T}{2}\right)\right),
\end{equation*}
\begin{equation*}
\left(q^{(1)}(0),q^{(2)}(0)\right)=\left(q^{(1)}\left(\frac{T}{2}\right),-q^{(2)}\left(\frac{T}{2}\right)\right),
\end{equation*}
that is
\begin{equation*}
q(0)=q\left(\frac{T}{2}\right)=\theta.
\end{equation*}
This means the original point $\theta=(0,0)$ is one of the points of
{\em self-intersection}s of $q(t)$. From the expression in
(\ref{e4}), one has
\begin{equation*}
q_{i}(t)=q\left(t+(N-i)\frac{T}{N}\right)=q\left(t-\frac{i}{N}T\right)
\end{equation*}
\begin{equation}\label{q3}
\quad\quad=\left(q^{(1)}(t-\frac{i}{N}T),q^{(2)}(t-\frac{i}{N}T)\right).
\end{equation}

Now, in the case of even bodies on Figure-Eight, i.e. $N=2k$, $k\in
\mathbb{Z}^{+}$, the following is from (\ref{q3}),
\begin{equation}\label{qi}
q_{i}(t)=\left(q^{(1)}\left(t-\frac{i}{2k}T\right),q^{(2)}\left(t-\frac{i}{2k}T\right)\right).
\end{equation}
For all $i\leq k$, we consider
\begin{eqnarray*}
q_{i+k}(t)&=&\left(q^{(1)}\left(t-\frac{i+k}{2k}T\right),q^{(2)}\left(t-\frac{i+k}{2k}T\right)\right)\\
&=&\left(q^{(1)}\left(t-\frac{i}{2k}T-\frac{T}{2}\right),q^{(2)}\left(t-\frac{i}{2k}T-\frac{T}{2}\right)\right).
\end{eqnarray*}
We notice this in (\ref{q1})
\begin{equation}\label{qik}
q_{i+k}(t)=\left(-q^{(1)}\left(t-\frac{i}{2k}T\right),q^{(2)}\left(t-\frac{i}{2k}T\right)\right).
\end{equation}
Together (\ref{qi}) with (\ref{qik}), it shows that when
$t=\frac{i}{2k}T$,
\begin{equation*}
q_{i}\left(\frac{i}{2k}T\right)=q_{i+k}\left(\frac{i}{2k}T\right)=\theta.
\end{equation*}
That is the $i$-th and $(i+k)$-th bodies will collide at the time
$t=\frac{i}{2k}T$.
\begin{proposition}
With the assumption of (\ref{e2}),(\ref{e4}) and (\ref{e5}), even
bodies on Figure-Eight orbit will collide on original point.
\end{proposition}

Since the constraint is more free, the shape of the possible orbit
is more complex. It is reasonable to make a further strict
assumption, a ``simple" orbits constraint, on Figure-Eight. We set
the further {\it Simple} Figure-Eight assumption:
\begin{equation}
\label{e7}
q(t_{1})=q(t_{2})\ \ \Rightarrow\ \ t_{i}\in\left\{0,\
\frac{T}{2}\right\},\quad i=1,2.
\end{equation}
It means that the choreography orbit of $q$ has one and only one
{\em self-intersection} on Figure-Eight.

Now we consider the case for the odd bodies on Figure-Eight. Let
$N=2k+1$, $k\in \mathbb{Z}^{+}$ and suppose that $q_{i}$ and $q_{j}$
collide at $t^{*}$, then
\begin{equation*}
q\left(t^{*}-\frac{i}{2k+1}T\right)=q\left(t^{*}-\frac{j}{2k+1}T\right).
\end{equation*}
According to assumption (\ref{e7}), one has
\begin{equation*}
t^{*}-\frac{i}{2k+1}T,\quad t^{*}-\frac{j}{2k+1}T\in\left\{0,
\frac{T}{2}\right\}.
\end{equation*}
Without losing generality, set
\begin{equation*}
t^{*}-\frac{i}{2k+1}T=0,\quad t^{*}-\frac{j}{2k+1}T=\frac{T}{2},
\end{equation*}
which means
\begin{equation*}
\frac{i-j}{2k+1}T=\frac{T}{2}.
\end{equation*}
Since $i, j, k\in \mathbb{Z}^{+}$, the expression above is
contradictory. Then we have
\begin{proposition}
With the assumption (\ref{e2}),(\ref{e4}), (\ref{e5}) and
(\ref{e7}), odd bodies will not collide on {\em simple} figure-eight
orbits.
\end{proposition}

It is obvious that assumption (\ref{e7}) can be generalized by
\begin{equation}\label{e71}
q(t_{1})=q(t_{2})\ \ \Rightarrow\ \ t_{i}\in I \quad i=1,2,
\end{equation}
where
\begin{equation*}
I=\left\{0\right\}\cup\left\{\frac{(2r_{1}-1)T}{2r_{2}}\ | \
r_{1},r_{2}\in\mathbb{Z}^{+}, \frac{2r_{1}-1}{2r_{2}}\leq 1\right\}.
\end{equation*}

\begin{theorem}\label{th1}
There is a dense subset in $[0,T]$, such that odd bodies on
Figure-Eight without collision in it.
\end{theorem}
{\bf Proof:} One can easily verify that $I$ is dense in $[0,T]$,
i.e., for $\forall t_{0}\in[0,T]$ and $\forall\epsilon$, there
exists $t^{*}\in I$, such that $ |t_{0}-t^{*}|\leq\epsilon$. Then
$I$ is at least one of the dense subset of $[0,T]$ we need. This
completes the proof.
\begin{remark}
If we could make a further generalization for assumption (\ref{e71})
by
\begin{equation*}
q(t_{1})=q(t_{2})\ \ \Rightarrow\ \ t_{i}\in [0,T],\quad i=1,2,
\end{equation*}
we would have gotten a {\em perfect} result. But we can not get the
generalization now. It means the Figure-Eight perhaps has some
self-intersection points on time $t^{*}$, which is not in $I$, e.g.
$t^{*}$ is an irrational point.
\end{remark}
\begin{remark}
It is obvious, when we record $t\in\mathbb{R}/T\mathbb{Z}$, a loop
space, we may have the same result by modifying
\begin{equation*}
I=\left\{0\right\}\cup\left\{\frac{(2r_{1}-1)T}{2r_{2}}\ | \
r_{1},r_{2}\in\mathbb{Z}^{+}\right\}.
\end{equation*}
\end{remark}

\begin{figure}[h]
    \centering
    \includegraphics[width=4in]{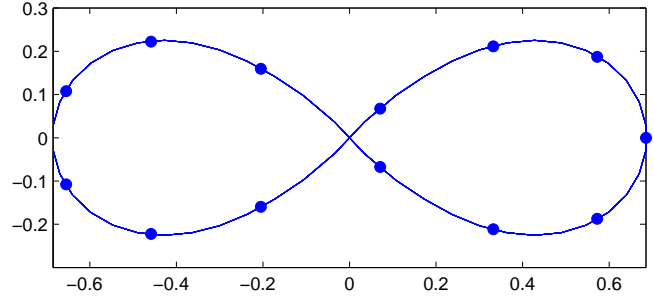}
    \caption{13 bodies on Figure-Eight.}\label{Eight13}
\end{figure}

{\bf Figure \ref{Eight13}} in Section \ref{se7} is one of our
numerical examples for odd bodies on Figure-Eight.

\section{Symmetries arising from Figure-Eight}\label{se4}

\subsection{From simple to double choreographies}\label{se41}

We consider some two choreographies orbits, the Double-Eight orbits,
which arises from Figure-Eight. With this idea, one can also figure
out other more choreographies orbits. We describe the following
Double-Eight orbits as an example for N-body problem.

Based on the aforementioned results, we set $N_{1}\geq 3$ and
$N_{2}\geq 3$ are all odd in $\mathbb{Z}^{+}$, and the number of
bodies is $N=N_{1}+N_{2}$. Then we make a two-choreographies
assumption: one Figure-Eight reclines on $x$-axis, the other on
$y$-axis. That is, one Figure-Eight $q(t)$ satisfies
\begin{equation}\label{D1}
q_{i}(t)=q\left(t+(N_{1}-i)\frac{T}{N_{1}}+\alpha_{1}\right),\quad
i=1,\cdots,N_{1}
\end{equation}
\begin{equation}\label{D2}
\sum^{N_{1}}_{i=1}q_{i}(t)=0,
\end{equation}
\begin{equation}\label{D3}
q(\sigma(t))=\sigma(q(t)),\quad q(\tau(t))=\tau(q(t)),
\end{equation}
the other Figure-Eight $Q(t)$ satisfies
\begin{equation}
\label{D5}
Q_{i}(t)=Q\left(t+(N_{2}-i)\frac{T}{N_{2}}+\alpha_{2}\right),\quad
i=1,\cdots,N_{2},
\end{equation}
\begin{equation}
\label{D6} \sum^{N_{2}}_{i=1}Q_{i}(t)=0,
\end{equation}
\begin{equation}\label{Q3}
Q(\sigma(t))=\tau(Q(t)),\quad Q(\tau(t))=\sigma(Q(t)),
\end{equation}
where $\alpha_{1}$ and $\alpha_{2}$ can be regarded as the {\em
phase difference}s after equivariant shifts. In order to guarantee
that any two bodies in two Figure-Eight separately will not collide
on the origin point, specially, we set $\alpha_{1}=0$,
$\alpha_{2}=\frac{T}{K}$, where $K\in\mathbb{Z}^{+}$ is prime to
$N_{1}$.

\begin{figure}[h]
    \begin{minipage}[t]{0.3\linewidth}
    \includegraphics[width=2in]{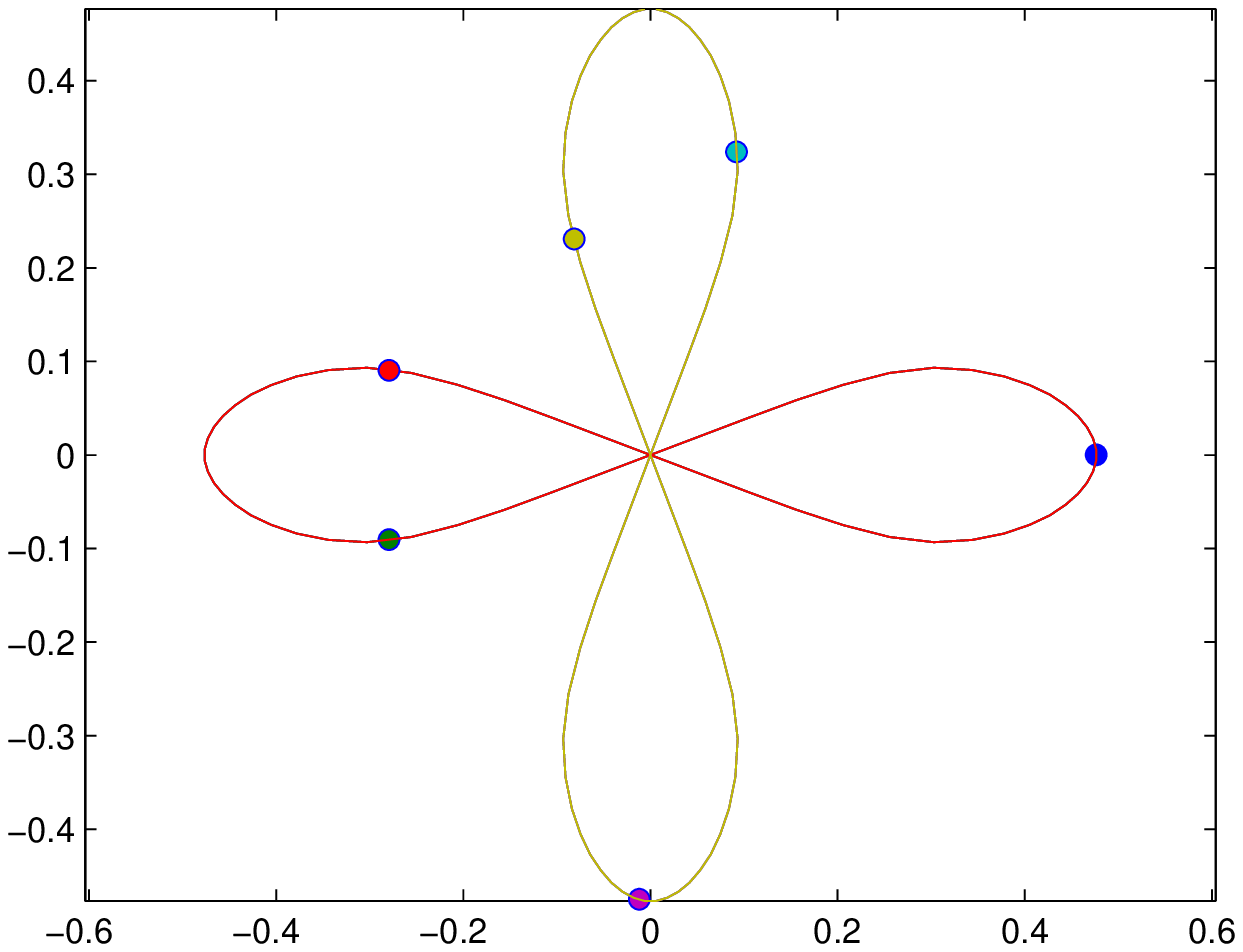}
    \caption{m=[1,1,1;1,1,1];\newline $\alpha_{1}=0, \alpha_{2}=\frac{T}{6}$.}\label{F3_3E01}
    \end{minipage}
    \hspace{0.3cm}
    \begin{minipage}[t]{0.3\linewidth}
    \includegraphics[width=2in]{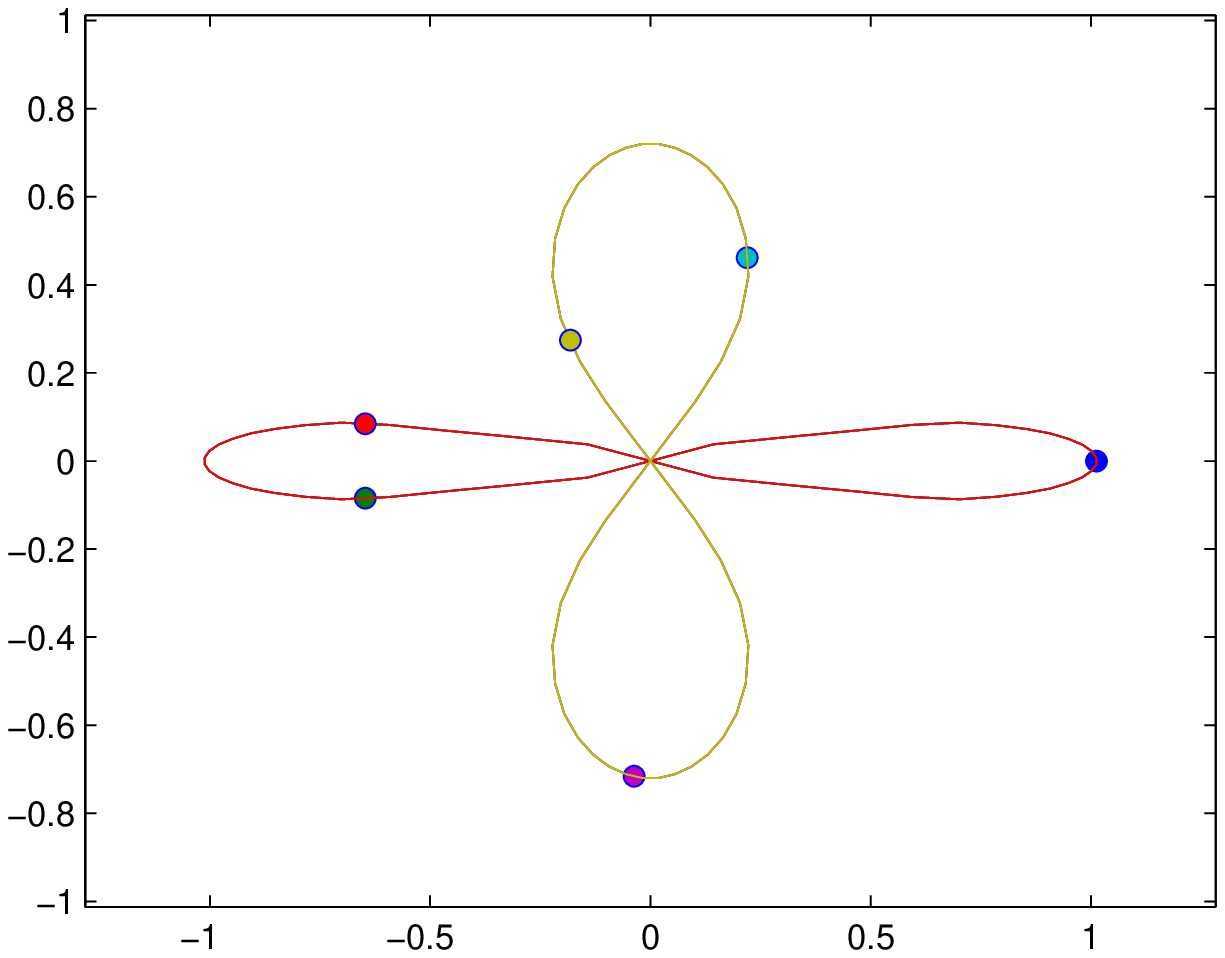}
    \caption{m=[1,1,1;10,10,10];\newline $\alpha_{1}=0, \alpha_{2}=\frac{T}{6}$.}\label{F3_3N01}
    \end{minipage}
    \hspace{0.3cm}
    \begin{minipage}[t]{0.3\linewidth}
    \includegraphics[width=2in]{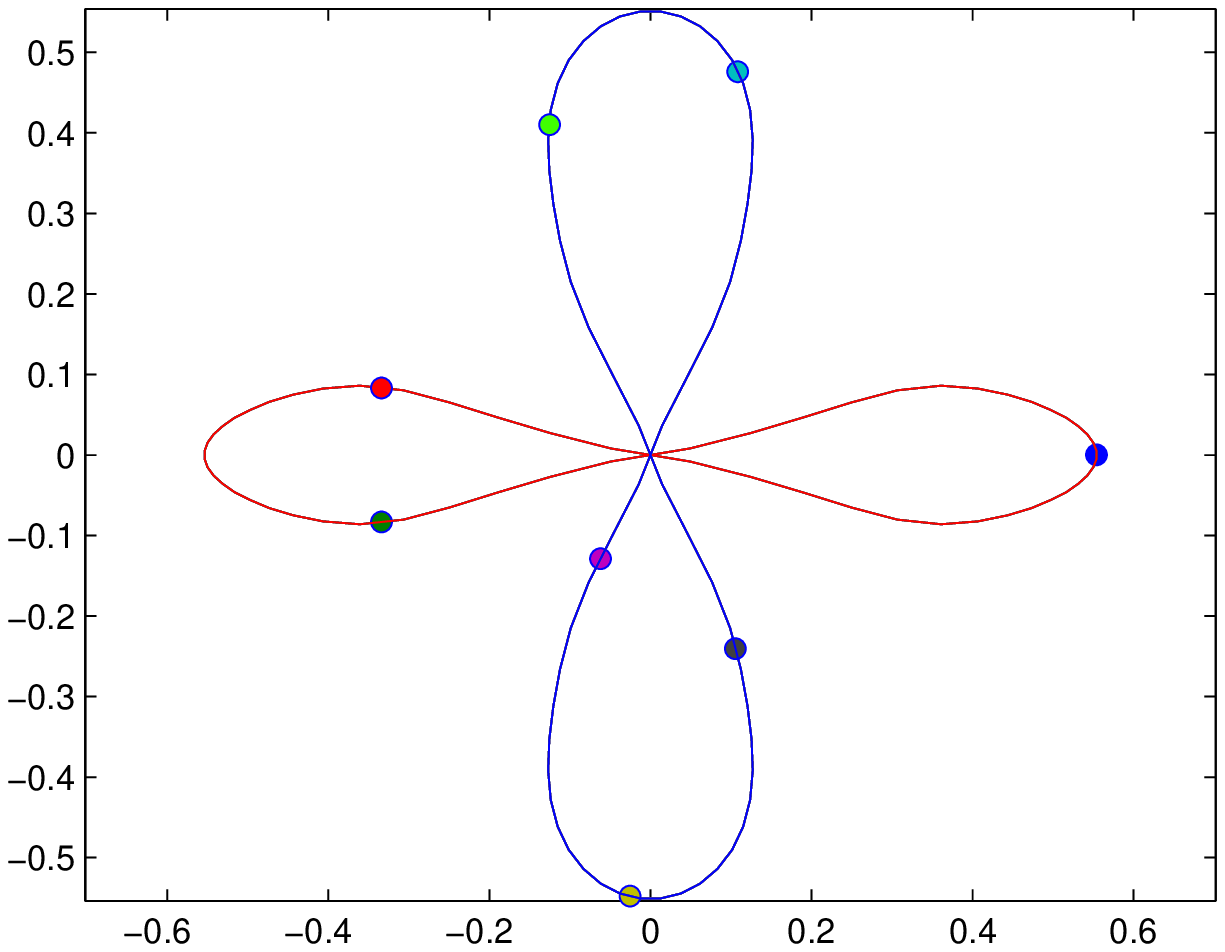}
    \caption{m=[1,1,1;1,1,1,1,1];\newline $\alpha_{1}=0, \alpha_{2}=\frac{T}{6}$.}\label{F3_5E01}
    \end{minipage}
\end{figure}

\begin{figure}[h]
    \begin{minipage}[t]{0.3\linewidth}
    \includegraphics[width=2in]{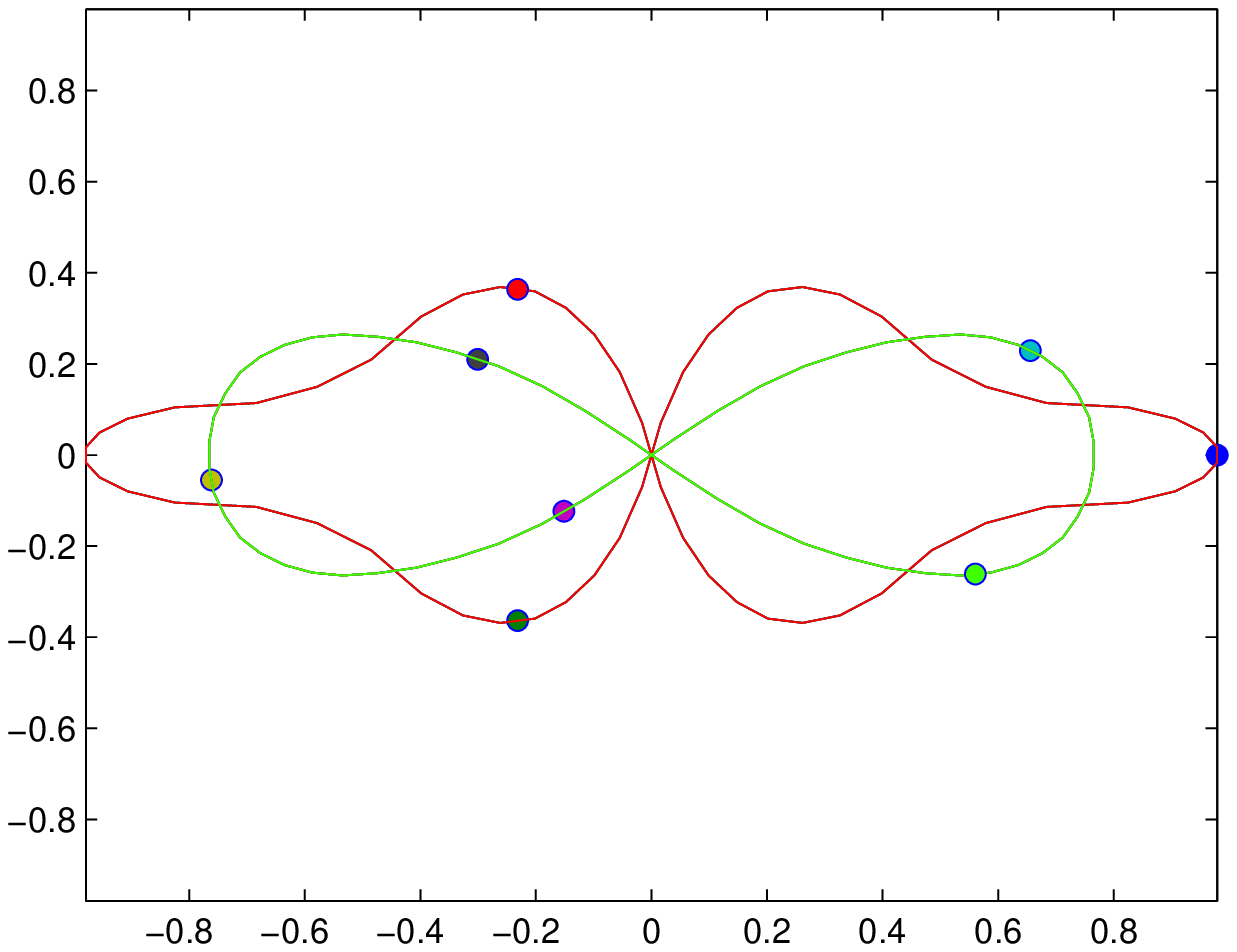}
    \caption{m=[1,1,1;5,5,5,5,5];\newline $\alpha_{1}=0, \alpha_{2}=\frac{T}{10}$.}\label{L3_5N04}
    \end{minipage}
    \hspace{0.3cm}
    \begin{minipage}[t]{0.3\linewidth}
    \includegraphics[width=2in]{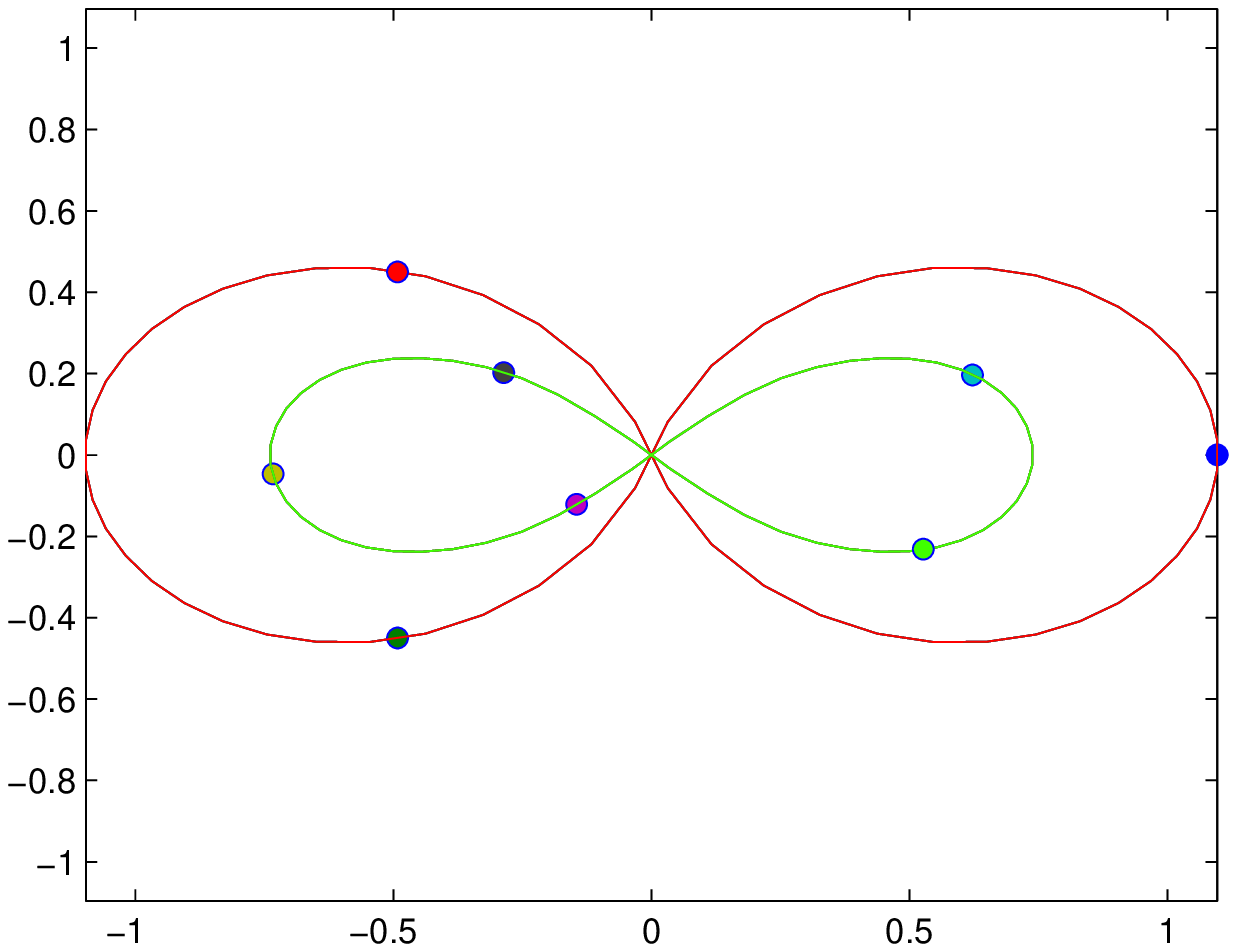}
    \caption{m=[1,1,1;5,5,5,5,5];\newline $\alpha_{1}=0, \alpha_{2}=\frac{T}{10}$.}\label{L3_5N07}
    \end{minipage}
    \hspace{0.3cm}
    \begin{minipage}[t]{0.3\linewidth}
    \includegraphics[width=2in]{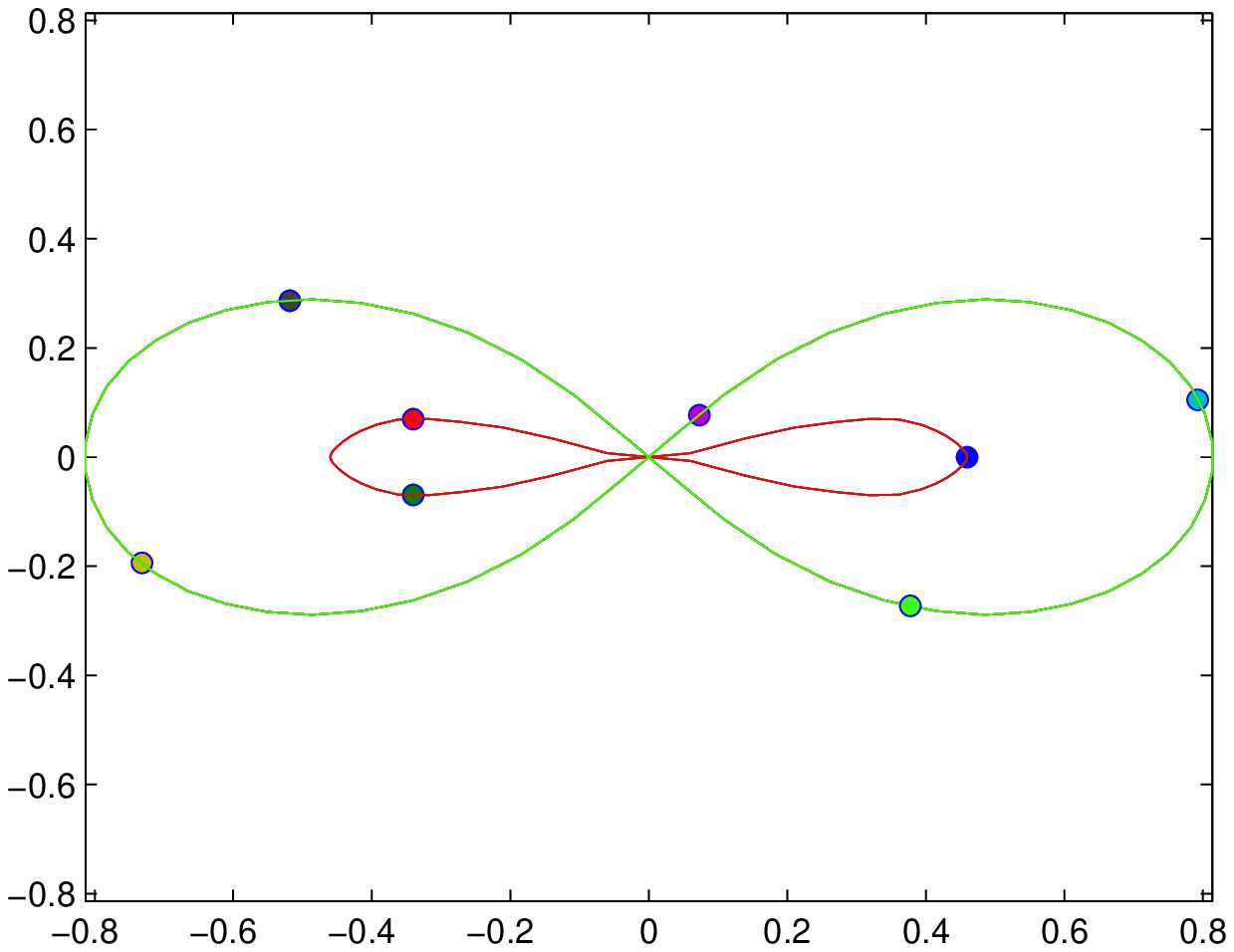}
    \caption{m=[1,1,1;5,5,5,5,5];\newline $\alpha_{1}=0, \alpha_{2}=\frac{T}{20}$.}\label{L3_5N05}
    \end{minipage}
\end{figure}

When we set $Q$ not satisfy (\ref{Q3}) but (\ref{D3}), we will get
the orbits with the constraint such that two Figure-Eights recline
on the same $x$-axis (or $y$-axis). The {\bf Figures
\ref{F3_3E01}}-{\bf \ref{L3_5N05}} in Section \ref{se7}, are the
numerical examples for Double-Eight orbits.

\subsection{Mixed-symmetry constraints in
$\mathbb{R}^{3}$}\label{se5} 
\begin{figure}[!ht]
    \begin{minipage}[t]{0.45\linewidth}
    \includegraphics[width=3in]{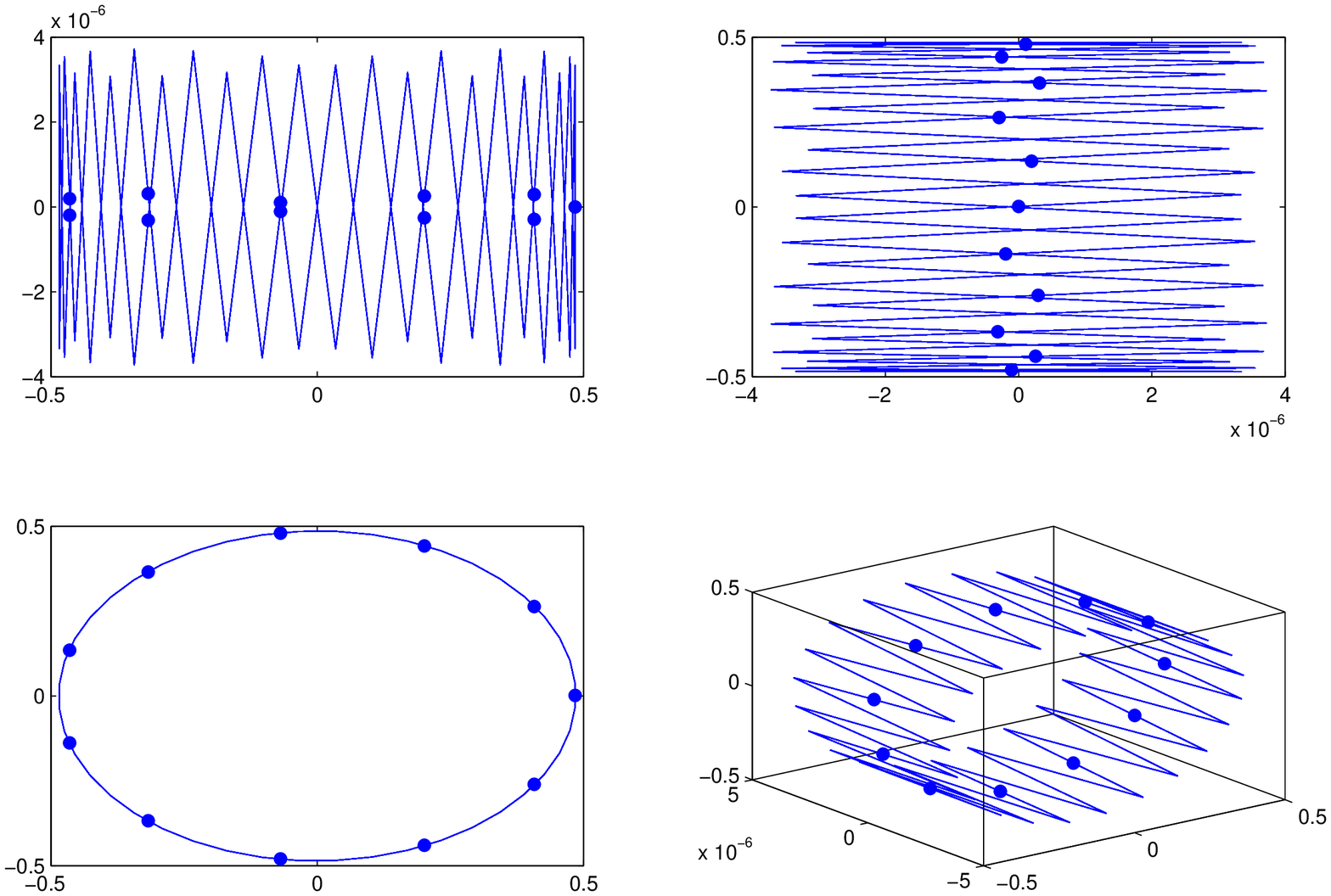}
    \caption{11 bodies with symmetry 1, \newline $m=[1,1,1,1,1,1,1,1,1,1,1];$}\label{R3X11E01}
    \end{minipage}
    \hspace{0.5cm}
    \begin{minipage}[t]{0.45\linewidth}
    \includegraphics[width=3in]{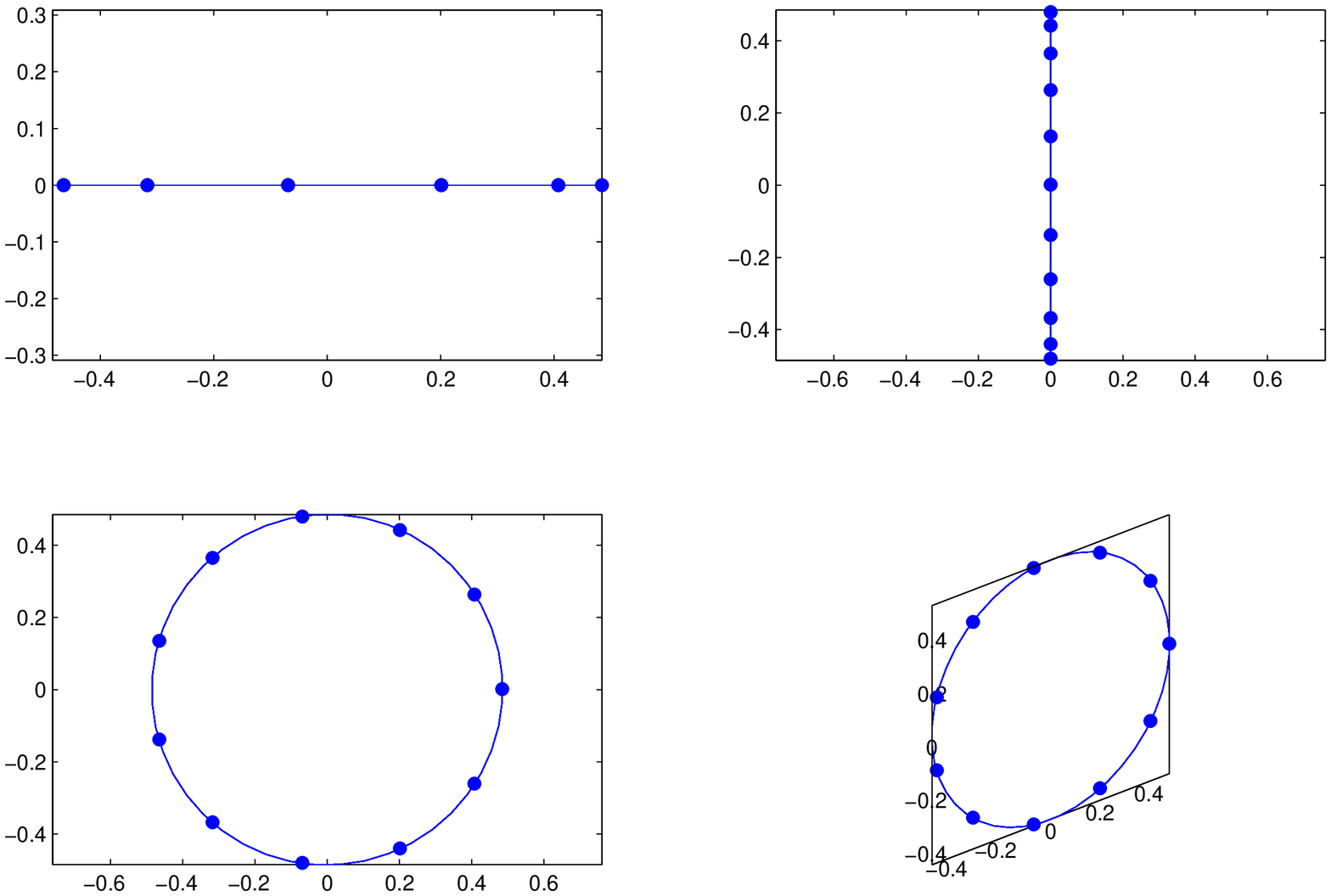}
    \caption{The same orbit as in Figure:\ref{R3X11E01} with different scale.}\label{Lagrange}
    \end{minipage}
\end{figure}

\begin{figure}[!ht]
    \begin{minipage}[t]{0.45\linewidth}
    \includegraphics[width=3in]{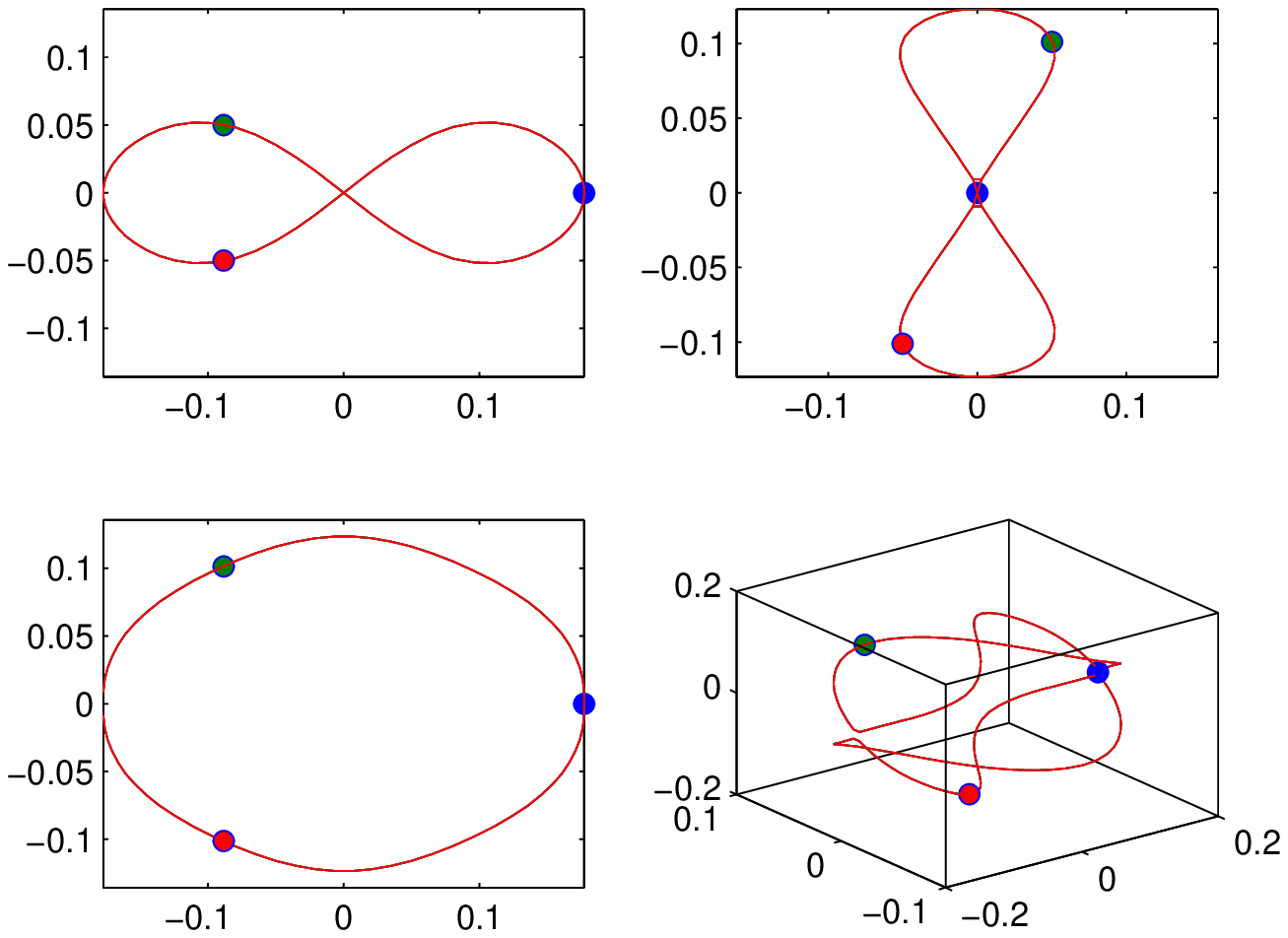}
    \caption{3 bodies with symmetry 2:\ \ \ \ \newline $m=[1,1,1].$}\label{M88N03E01}
    \end{minipage}
    \hspace{0.5cm}
    \begin{minipage}[t]{0.45\linewidth}
    \includegraphics[width=3in]{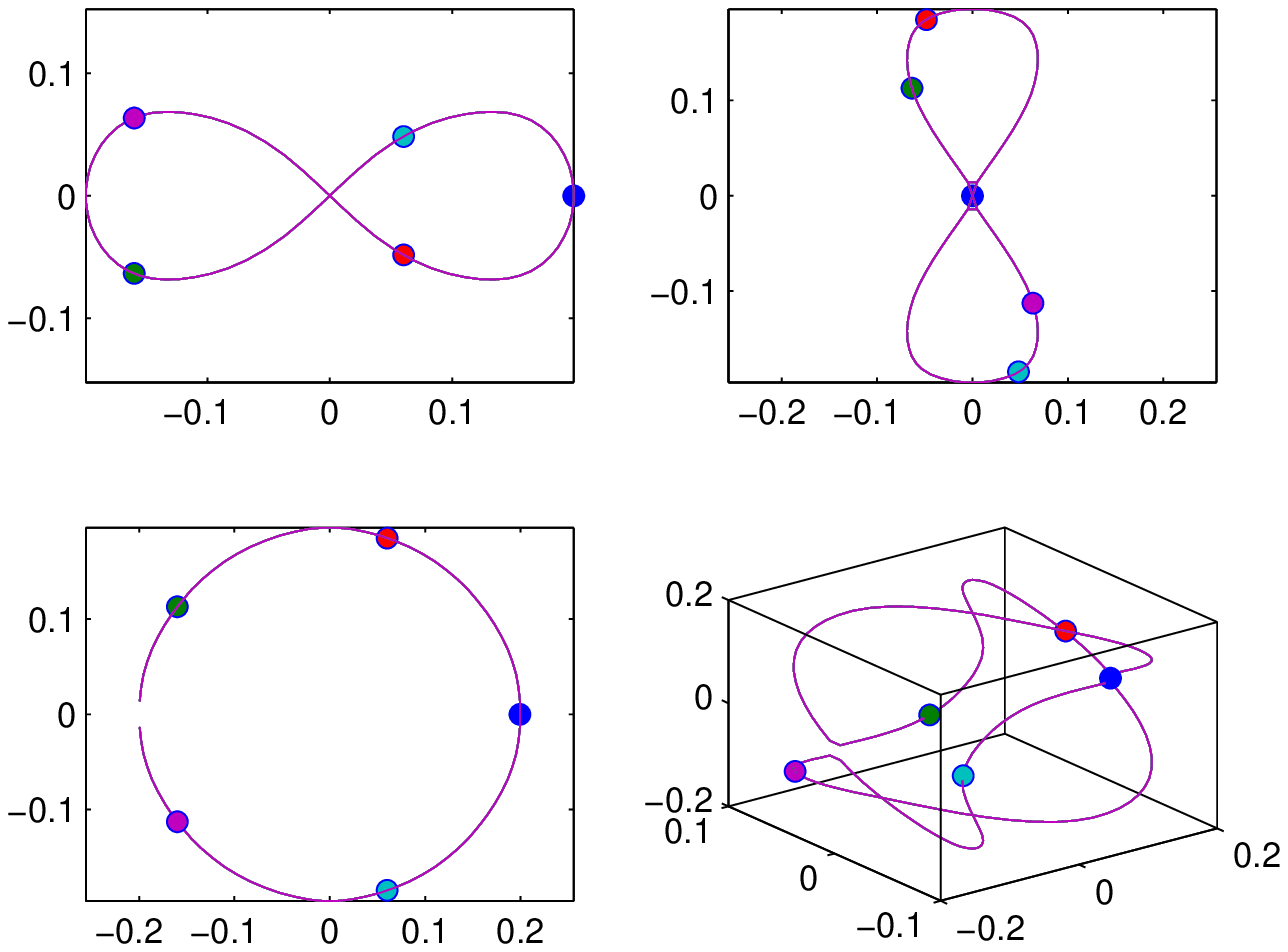}
    \caption{5 bodies with symmetry 2:\newline $m=[1,1,1,1,1].$}\label{M88N05E01}
    \end{minipage}
\end{figure}

\begin{figure}[!ht]
    \begin{minipage}[t]{0.45\linewidth}
    \includegraphics[width=3in]{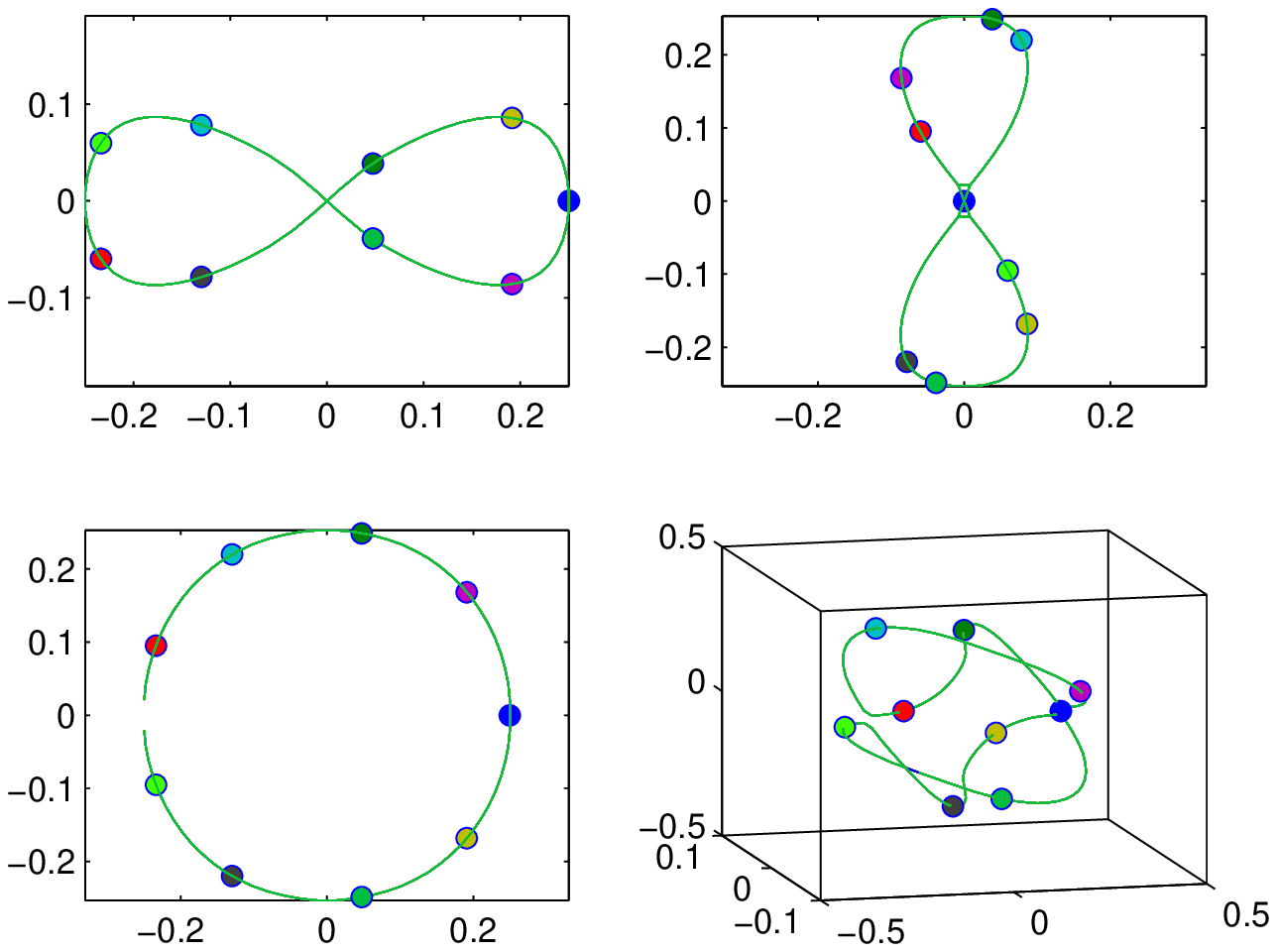}
    \caption{9 bodies with symmetry 2:\newline $m=[1,1,1,1,1,1,1,1,1].$}\label{M88N09E01}
    \end{minipage}
    \hspace{0.5cm}
    \begin{minipage}[t]{0.45\linewidth}
    \includegraphics[width=3in]{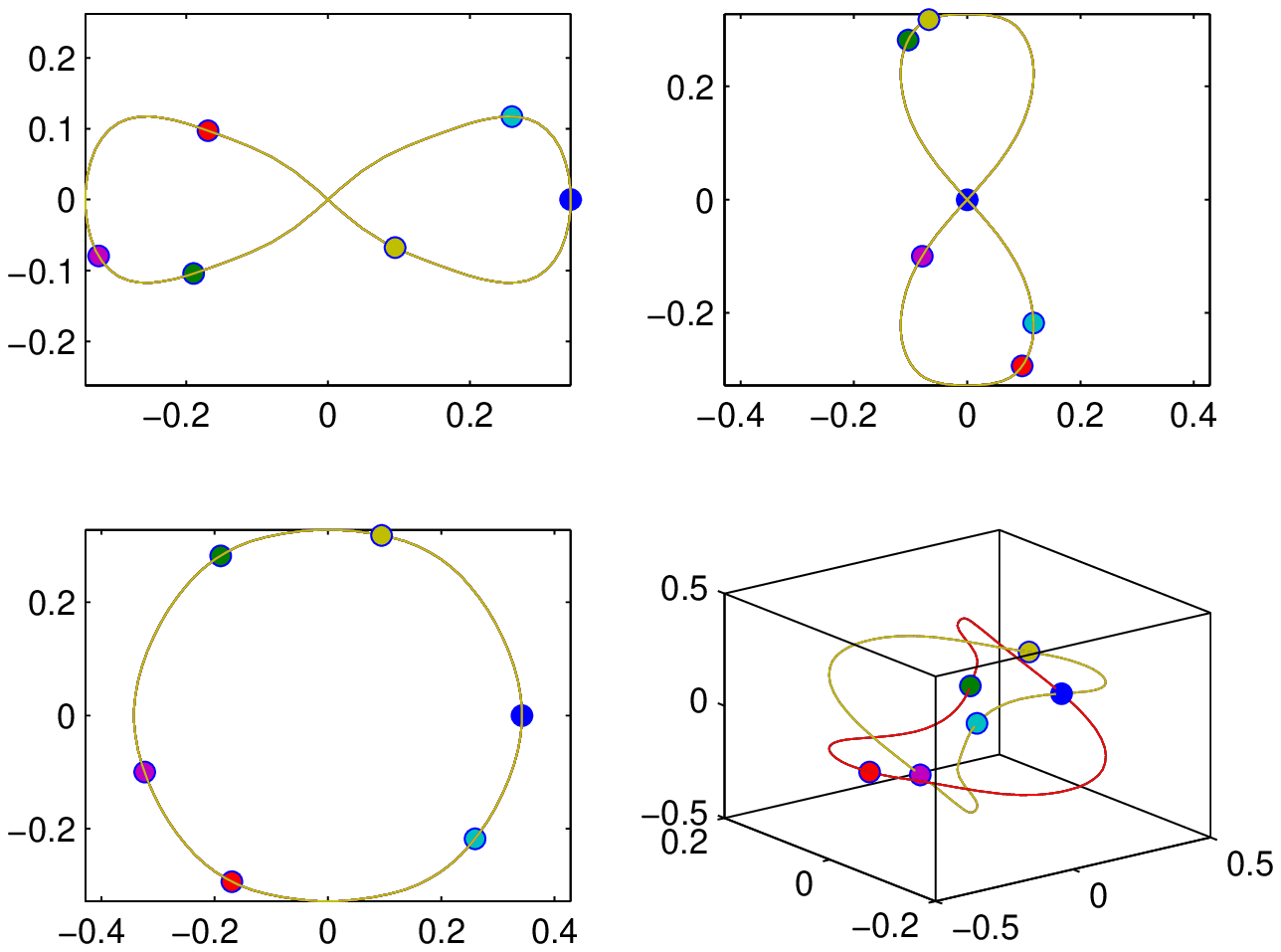}
    \caption{6 bodies with symmetry 3:\newline $m=[1,1,1;1,1,1], \alpha_{1}=0, \alpha_{2}=\frac{T}{8}$.}\label{M8ccl06E01}
    \end{minipage}
\end{figure}

\begin{figure}[!ht]
    \begin{minipage}[t]{0.45\linewidth}
    \includegraphics[width=3in]{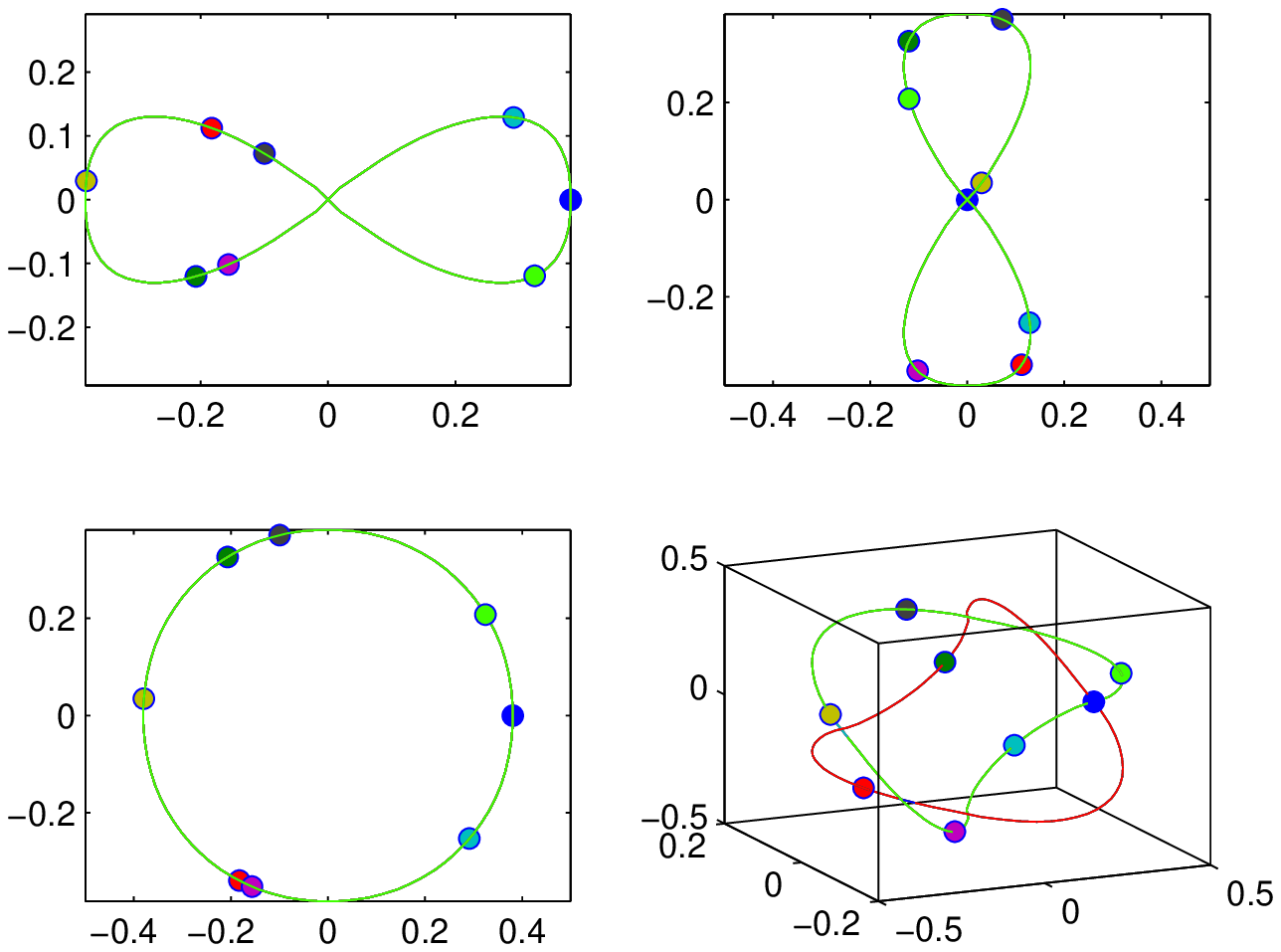}
   \caption{8 bodies with symmetry 3:\newline $m=[1,1,1;1,1,1,1,1], \alpha_{1}=0, \alpha_{2}=\frac{T}{8}$}\label{M8ccl08E01}
    \end{minipage}
    \hspace{0.5cm}
    \begin{minipage}[t]{0.45\linewidth}
    \includegraphics[width=3in]{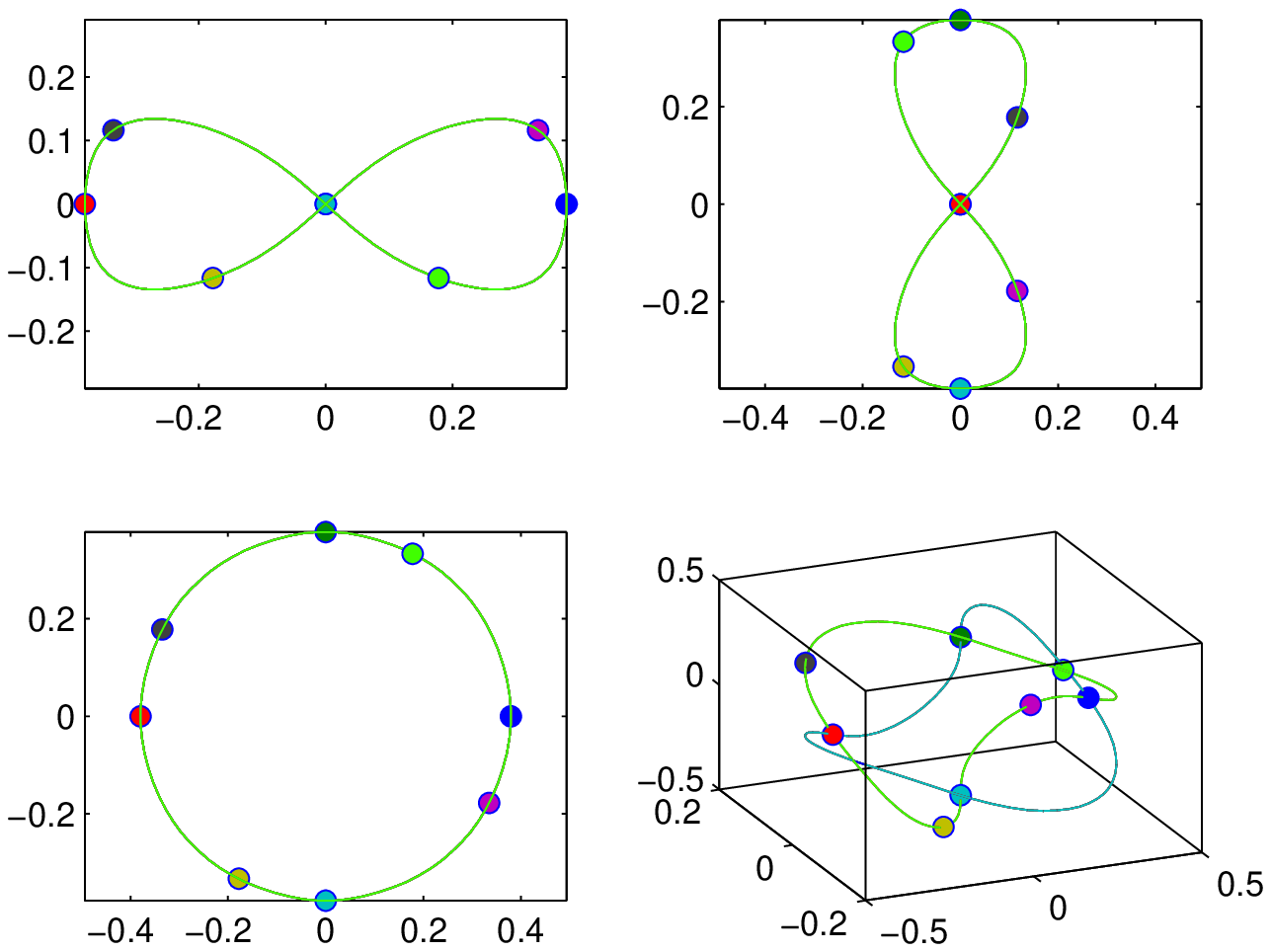}
    \caption{8 bodies with symmetry 3:\newline $m=[1,1,1,1;1,1,1,1], \alpha_{1}=0, \alpha_{2}=\frac{T}{12}$}\label{M8ccl08E02}
    \end{minipage}
\end{figure}

To construct some Mixed-symmetric
constraints on simple choreography or two choreographies for
$q=(q^{(1)},q^{(2)},q^{(3)})\in\mathbb{R}^{3}$, we set $q^{(1)}$ and
$q^{(2)}$ satisfy (\ref{D1})-(\ref{D3}) on $x-y$ plane, and in
$z$-coordinate, $q^{(3)}$, satisfy one of the following constraints:
\begin{itemize}
\item Symmetry 1 (simple choreography orbit)
\begin{equation*}
q^{(3)}(t)=-q^{(3)}\left(t+\frac{T}{2}\right),\quad
t\in\left[0,T\right].
\end{equation*}

\item Symmetry 2 (simple choreography orbit)
\begin{equation}\label{con1}
q^{(3)}(t)=q^{(3)}\left(\frac{T}{2}-t\right),\quad
t\in\left[\frac{T}{4},\frac{T}{2}\right],
\end{equation}
and
\begin{equation}\label{con2}
q^{(3)}(t)=-q^{(3)}\left(t+\frac{T}{2}\right),\quad
t\in\left[0,T\right].
\end{equation}

\item Symmetry 3 (two choreographies orbits)

Let $q=(q^{(1)},q^{(2)},q^{(3)})$ be a simple choreography and
$Q=(Q^{(1)},Q^{(2)},Q^{(3)})$ another one, we set
\begin{enumerate}
\item
$q$ and $Q$ locate on the same Figure-Eight on $x-y$ plane;
\item
$q^{(3)}$ satisfy (\ref{con1}) and (\ref{con2});
\item
$Q^{(3)}=-q^{(3)}$.
\end{enumerate}

\end{itemize}
See {\bf Figures \ref{R3X11E01}-\ref{M8ccl08E02}} for numerical
examples.

\subsection{Return to Figure-Eight} From the discussion in Section
\ref{se3} and Section \ref{se4}, we can choose many $\alpha_{i}$,
$i=1,\cdots,M$, $M\in\mathbb{Z}^{+}$, such that $M$ groups of odd
bodies can be arranged on Figure-Eight without any collision in a
dense subset of $[0,T]$.
\begin{theorem}\label{th3}
For all $N\geq 3$, $N\in\mathbb{Z}^{+}$, $N$ bodies can be arranged
on Figure-Eight without any collision in a dense subset of $[0,T]$.
\end{theorem}
{\bf Proof:} If $N\geq 3$ is odd, by Theorem \ref{th1}, the proof is
as follows. If $N$ is even, it can be decomposed by two odds $N_{1}$
and $N_{2}$, such that $N=N_{1}+N_{2}$. Then we set
\begin{equation*}
q_{i}(t)=q\left(t+(N_{1}-i)\frac{T}{N_{1}}+\alpha_{1}\right),\quad
i=1,\cdots,N_{1}
\end{equation*}
\begin{equation*}
\label{D5}
Q_{i}(t)=q\left(t+(N_{2}-i)\frac{T}{N_{2}}+\alpha_{2}\right),\quad
i=1,\cdots,N_{2},
\end{equation*}
\begin{equation*}
\sum^{N_{1}}_{i=1}q_{i}(t)+\sum^{N_{2}}_{i=1}Q_{i}(t)=0,
\end{equation*}
\begin{equation*}
q(\sigma(t))=\sigma(q(t)),\quad q(\tau(t))=\tau(q(t)),
\end{equation*}
\begin{equation*}
Q(\sigma(t))=\sigma(Q(t)),\quad Q(\tau(t))=\tau(Q(t)),
\end{equation*}
\begin{equation*}
\alpha_{1}=0,\quad\alpha_{2}=\frac{T}{K},
\end{equation*}
where $K\in\mathbb{Z}^{+}$ is prime to $N_{1}$. Through the
discussion in the counterpart of Section \ref{se4}, one can easily
verify that the theorem follows.

{\bf Figures \ref{Even06E01}-\ref{Even12E01}} are numerical
examples for Theorem \ref{th3}.

\section{Existence of symmetric orbits}\label{se6}
In this section, we mainly consider the existence of generalized
solutions and non-collision solutions for Mixed-symmetry and
Double-Eight.

At first, we assemble some well-known facts from \cite{Go1},
\cite{Go2}, \cite{MW} and \cite{Z4} about Sobolev space
$\mathbb{H}^{1}(\mathbb{R}/T\mathbb{Z}, (\mathbb{R}^{d})^{N})$ and
the ordinary space $\mathbb{C}^{2}([0,T], (\mathbb{R}^{d})^{N})$.

\begin{itemize}
\item
Let $\xi_{1}, \xi_{2}\in\mathbb{C}^{2}([0,T],
(\mathbb{R}^{d})^{N})$, then the inner product can be defined by
\begin{equation*}
\langle\xi_{1},\xi_{2}\rangle_{0}=\int_{0}^{T}\langle\xi_{1}(t),\xi_{2}(t)\rangle
dt,
\end{equation*}
and the corresponding norm is $\parallel \xi\parallel _{0}$.

\item
Let $\xi_{1}, \xi_{2}\in\mathbb{H}^{1}(\mathbb{R}/T\mathbb{Z},
(\mathbb{R}^{d})^{N})$, then the inner product can be defined by
\begin{equation*}
\langle\xi_{1},\xi_{2}\rangle_{1}=\langle\xi_{1},\xi_{2}\rangle_{0}+\langle\dot{\xi}_{1},\dot{\xi}_{2}\rangle_{0},
\end{equation*}
and the corresponding norm is
$\parallel\xi\parallel_{1}=\sqrt{\parallel\xi\parallel_{0}^{2}+\parallel\dot{\xi}\parallel_{0}^{2}}$.

\item By Sobolev imbedding theorems, weak convergence in
$\parallel\cdot\parallel_{1}$ implies uniform convergence, i.e., the
weak $\mathbb{H}^{1}$ topology is stronger than the $\mathbb{C}^{0}$
topology.

\item
Let $J$ be a functional from $\mathbb{H}^{1}$ to $\mathbb{R}$. A
sequence $\left\{^{k}q\right\}_{k=1}^{\infty}$ is called a {\em
minimizing sequence}, if
\begin{equation*}
J(^{k}q)\rightarrow\inf J, \quad\mbox{ whenever}\quad
k\rightarrow\infty.
\end{equation*}

\item
A functional $J:\mathbb{H}^{1}\rightarrow\mathbb{R}$ is {\em weakly
lower semi-continuous} if
\begin{equation*}
^{k}q\rightharpoonup q \quad \Rightarrow \quad \underline{\lim}
J(^{k}q)\geq J(q).
\end{equation*}

\item
A weakly lower semi-continuous functional $J$ on $\mathbb{H}^{1}$
has a minimum in $\mathbb{H}^{1}$, if $J$ has a bounded minimizing
sequence.
\end{itemize}

Since it is difficult to exclude the collision without any more
additional constraints for the orbits, we introduce the following
two definitions, both of which partially exclude the collision in
some sense.

We denote
\begin{equation*}
\Delta =\left\{t\in[0,T] \ | \ q_{i}(t)=q_{j}(t),\ \mbox{for some}\
i \neq j \right\}.
\end{equation*}

\begin{definition}
We call $q=(q_{1},\cdots,q_{N})\in \mathbb{C}^{2}([0,T],
(\mathbb{R}^{d})^{N})$ a {\em non-collision solution} of (\ref{N})
and (\ref{P}), if $\Delta=\phi$.
\end{definition}
\begin{definition}\label{defi2}
We call $q=(q_{1},\cdots,q_{N})\in
\mathbb{H}^{1}(\mathbb{R}/T\mathbb{Z}, (\mathbb{R}^{d})^{N})$ a {\em
generalized solution} of (\ref{N}) and (\ref{P}), if

a) $meas(\Delta)=0$;

b) $q$ satisfy (\ref{N});

c)
\begin{equation}\label{def2}
 \frac{1}{2}\sum_{i=1}^{N}m_{i}\left|\dot{q}_{i}(t)\right|^{2}-U(q(t))\equiv
C,
\end{equation}
where $t\in[0,T]\backslash \Delta$ and $C$ is a constant.
\end{definition}

In order to overcome the singularity of $U$, the following condition
is usually used as an efficient method in considering the existence
of minimum of functional $J$.
\begin{definition} Let some $V_{ij}\in \mathbb{C}^{1}(R^{d} \backslash \{0 \} , R )$,
$V_{ij}\rightarrow +\infty$ ($\xi \rightarrow 0$) and
\begin{equation*}
-U_{ij}(\xi)\geq |\nabla V_{ij}(\xi)|^{2}, \qquad \forall \xi \in
R^{d} \backslash \{0 \},
\end{equation*}
where $|\xi|$ is small. Then we say $U_{ij}$ satisfies {\em Strong
Force {\rm (SF)} condition}.
\end{definition}

For the individuality of the Mixed-symmetric (the same as
Double-Eight) solutions, we have to consider them respectively. But
in the following theorem, we describe our result at one time for all
the symmetric constraints, which have been mentioned therein.

\begin{theorem}\label{th2}
There are infinitely many generalized solutions for
(\ref{N})-(\ref{P}) with anyone of those symmetric constraints.
Furthermore, if $U$ satisfies (SF) condition, there will be
infinitely many non-collision solutions.
\end{theorem}

The proof of theorem \ref{th2} is extremely along the line from
Gordon \cite{{Go1},{Go2}} and Zelati \cite{Z4}. The key points
locate on the estimation of $|q(t)|$ and the construction of those
nested intervals to cover $[0,T]\backslash\Delta$. In the ongoing
proof, the existence of periodic solutions with symmetry 1 is proved
as a candidate of those Mixed-symmetries and Double-Eight. \ \
\newline{\bf Proof:}
We set
\begin{equation*}
\Omega=\left\{q \ | \
q=\left(q_{1},\cdots,q_{N}\right)\in(\mathbb{R}^{d})^{N}, q_{i}\neq
q_{j}, \mbox{for all }i\neq j\right\},
\end{equation*}
\begin{equation*}
\Lambda=\left\{q \ | \
q\in\mathbb{H}^{1}(\mathbb{R}/T\mathbb{Z},\Omega)\right\},
\end{equation*}
\begin{equation*}
\Lambda_{0}=\left\{q \ | \ q\in\Lambda,\quad q \ \mbox{satisfies
symmetry 1}\right\}.
\end{equation*}
It's easy to verify that the critical points of $J$ on $\Lambda$ are
non-collision solutions of (\ref{N})-(\ref{P}). And one can easily
check that the critical points of $J|_{\Lambda_{0}}$ are actually
those of $J$ on $\Lambda$.

{\bf Step1:} If $U_{ij}$ does not satisfy (SF) condition, we can
modify it by the following formula
\begin{equation*}
U^{\delta}_{ij}=U_{ij}-\frac{\phi(|q_{j}-q_{i}|)}{|q_{j}-q_{i}|^{2}},
\end{equation*}
where $\delta$ is small and
\begin{displaymath}
\phi = \left\{
    \begin{array}{ll}
    0, & \quad\mbox{when} \ |q_{j}-q_{i}|\leq \displaystyle\frac{\delta}{2},\\
    1, & \quad\mbox{when} \ |q_{j}-q_{i}|\geq \delta.
    \end{array} \right.
\end{displaymath}
Then the modified $U^{\delta}$ satisfies (SF) condition, and the
corresponding Lagrangian functional is
\begin{equation*}
J^{\delta}(q)=\frac{1}{2}\sum_{i=1}^{N}m_{i}\int_{0}^{T}|\dot{q}_{i}|^{2}dt-\int_{0}^{T}U^{\delta}(q)dt,
\end{equation*}
where
\begin{equation*}
U^{\delta}(q)=\frac{1}{2}\sum_{1\leq i\neq j\leq
N}U_{ij}^{\delta}(q_{j}-q_{i}).
\end{equation*}

{\bf Step2:} In this step, we will show the existence of a minimum
of $J$. Let the infimum of $J^{\delta}$ be
\begin{equation*}
c_{\delta}=\inf \{ \ J^{\delta}(q) \ | \ q\in \Lambda_{0} \ \}.
\end{equation*}
Consider a minimizing sequence $^{k}q\in \Lambda_{0}$ such that
$J^{\delta}(^{k}q)\rightarrow c_{\delta}$ when $k\rightarrow
+\infty$. Then, when $k$ is large enough, we have
\begin{equation*}\label{Jd}
J^{\delta}(^{k}q(t))\leq c_{\delta}+1.
\end{equation*}
Hence for one term of $J^{\delta}$, we can deduce that
\begin{equation*}
\int_{0}^{T}\left|^{k}\dot{q}_{i}\right|^{2}dt\leq
\frac{2(c_{\delta}+1)}{\tilde{m}}\triangleq C, \quad\forall i,
\end{equation*}
where $C$ is a constance and $\tilde{m}=\min\{m_{1},\cdots,m_{N}\}$.

We notice that from (\ref{q1})-(\ref{q2}), one can deduce the period
of $q^{(2)}$ is $\frac{T}{2}$, and $q(t)$ is symmetrical about
$q^{(2)}$-axis, then one has
\begin{equation}
q^{(2)}(t)=-q^{(2)}\left(t+\frac{T}{4}\right).
\end{equation}
Hence we can estimate the norm of $q$,
\begin{eqnarray*}
\left|q(t)\right| & \leq &
\left|q^{(1)}(t)\right|+\left|q^{(2)}(t)\right|+\left|q^{(3)}(t)\right| \\
& =  & \frac{1}{2}\left|q^{(1)}(t)-q^{(1)}(t+\frac{T}{2})\right|
+\frac{1}{2}\left|q^{(2)}(t)-q^{(2)}(t+\frac{T}{4})\right|\\
&    &
+\frac{1}{2}\left|q^{(3)}(t)-q^{(3)}(t+\frac{T}{2})\right|\\
& =  &
\frac{1}{2}\left|\int_{t}^{t+\frac{T}{2}}\dot{q}^{(1)}(s)ds\right|
+\frac{1}{2}\left|\int_{t}^{t+\frac{T}{4}}\dot{q}^{(2)}(s)ds\right|
+\frac{1}{2}\left|\int_{t}^{t+\frac{T}{2}}\dot{q}^{(3)}(s)ds\right|\\
&\leq&
\frac{1}{2}\int_{t}^{t+\frac{T}{2}}\left|\dot{q}^{(1)}(s)\right|ds
+\frac{1}{2}\int_{t}^{t+\frac{T}{4}}\left|\dot{q}^{(2)}(s)\right|ds
+\frac{1}{2}\int_{t}^{t+\frac{T}{2}}\left|\dot{q}^{(3)}(s)\right|ds\\
&\leq&
\frac{1}{2}\left(\int_{t}^{t+\frac{T}{2}}\left|\dot{q}^{(1)}(s)\right|^{2}ds\right)^{\frac{1}{2}}
\left(\frac{T}{2}\right)^{\frac{1}{2}}\\
&    &
+\frac{1}{2}\left(\int_{t}^{t+\frac{T}{4}}\left|\dot{q}^{(2)}(s)\right|^{2}ds\right)^{\frac{1}{2}}
\left(\frac{T}{4}\right)^{\frac{1}{2}}\\
&    &
+\frac{1}{2}\left(\int_{t}^{t+\frac{T}{2}}\left|\dot{q}^{(3)}(s)\right|^{2}ds\right)^{\frac{1}{2}}
\left(\frac{T}{2}\right)^{\frac{1}{2}}\\
&\leq&
\frac{\sqrt{T}}{4}\left(\left(\int_{0}^{T}\left|\dot{q}^{(1)}(t)\right|^{2}dt\right)^{\frac{1}{2}}
+\left(\int_{0}^{T}\left|\dot{q}^{(2)}(t)\right|^{2}dt\right)^{\frac{1}{2}}\right.\\
&    &
\left.+\left(\int_{0}^{T}\left|\dot{q}^{(3)}(t)\right|^{2}dt\right)^{\frac{1}{2}}\right)\\
&\leq&
\frac{\sqrt{3T}}{4}\left(\int_{0}^{T}\left|\dot{q}^{(1)}(t)\right|^{2}dt
+\int_{0}^{T}\left|\dot{q}^{(2)}(t)\right|^{2}dt+\int_{0}^{T}\left|\dot{q}^{(3)}(t)\right|^{2}dt\right)^{\frac{1}{2}}\\
&=&
\frac{\sqrt{3T}}{4}\left(\int_{0}^{T}\left|\dot{q}(t)\right|^{2}dt\right)^{\frac{1}{2}}.
\end{eqnarray*}

Then for $\forall i$, $\forall k$ (sufficiently large), the norm for
one of the minimizing sequence is
\begin{equation*}
\parallel^{k}q\parallel_{1}=\sqrt{\parallel^{k}q\parallel_{0}^{2}+\parallel^{k}\dot{q}\parallel_{0}^{2}}\leq
\frac{1}{4}\sqrt{3TC(T+1)}.
\end{equation*}
This implies the existence of
$\bar{q}^{\delta}=(\bar{q}_{1}^{\delta},\cdots,\bar{q}_{N}^{\delta})$
with
$\bar{q}_{i}^{\delta}\in\mathbb{H}^{1}(\mathbb{R}/T\mathbb{Z},\mathbb{R}^{d})$,
for all $i=1,\cdots,N$ such that
\begin{equation*}
^{k}q_{i}\rightharpoonup\bar{q}_{i}^{\delta}, \quad\forall\ \
^{k}q_{i}\in\mathbb{H}^{1}(\mathbb{R}/T\mathbb{Z},\mathbb{R}^{d}),
\end{equation*}
and
\begin{equation*}
^{k}q_{i}\rightarrow\bar{q}_{i}^{\delta}, \quad\forall\ \
^{k}q_{i}\in\mathbb{C}^{0}(\mathbb{R}/T\mathbb{Z},\mathbb{R}^{d}).
\end{equation*}

On condition (SF), it follows that $J(^{k}q)\rightarrow +\infty$ for
every sequence $\{^{k}q\}$ such that $^{k}q\rightarrow\bar{q}$
weakly in $\mathbb{H}^{1}$ and strongly in $\mathbb{C}^{0}$ if
$\bar{q}\in\partial\Lambda$. This proves that
$\bar{q}^{\delta}\in\Lambda_{0}$ is a minimum for $J^{\delta}$ on
$\Lambda_{0}$. This minimum is a non-collision solution of
(\ref{N})-(\ref{P}). Then we have proved that (\ref{N})-(\ref{P})
has at least one non-collision solution if condition (SF) is
satisfied.

{\bf Step3:} This step is the same as the corresponding part of the
proof for Theorem1.1 in \cite{Z4}, so we only describe the sketch:
We can find a generalized solution $q^{*}$ of (\ref{N})-(\ref{P}) by
constructing compact nested interval cover of $[0,T]\backslash
\Delta$ to show $U^{\delta}\rightarrow U$ and
$q^{*}_{\delta}\rightarrow q^{*}$ when $\delta\rightarrow 0$. This
process shows that $q^{*}$ satisfies definition \ref{defi2}.

We can also prove $q_{T/z}$ is a $T$-periodic solution of
(\ref{N})-(\ref{P}) if $q$ is a $T$-periodic solution, where $z$ is
a positive integer, then we can find infinitely many $T$-periodic
solutions of (\ref{N})-(\ref{P}).  This completes the proof.

\section{Numerical examples}\label{se7}
In this section, we along the line in Theorem \ref{th2}, finding
the critical points of $J$, to list some of numerical examples for
those theoretic consideration. In the following figures, $m$ is
the masses, some of which are divided into two parts with a
semicolon to distinguish the masses of $N_{1}$ bodies and $N_{2}$
bodies, $\alpha_{1}$ and $\alpha_{2}$ as defined above. In those
figures with sub-figures, the top left sub-figure is the
projection of the orbits on $x$-$y$ plane, the top right is the
projection of the orbits on $y$-$z$ plane, the bottom left is the
projection of the orbits on $x$-$z$ plane, and the bottom right is
the orbits in space.

\begin{figure}[!ht]
    \begin{minipage}[t]{0.45\linewidth}
    \includegraphics[width=3in]{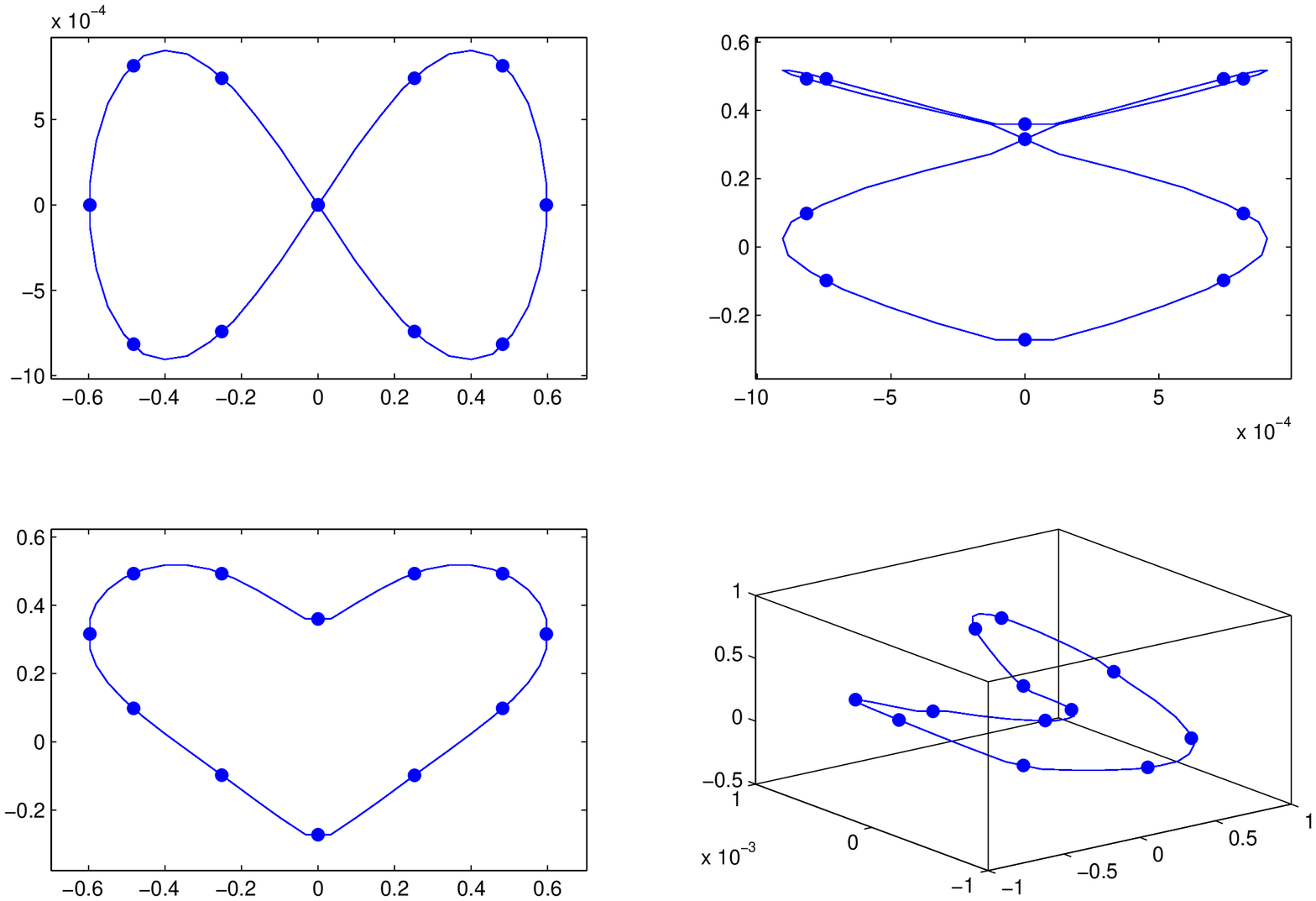}
    \end{minipage}
    \hspace{0.5cm}
    \begin{minipage}[t]{0.45\linewidth}
    \includegraphics[width=3in]{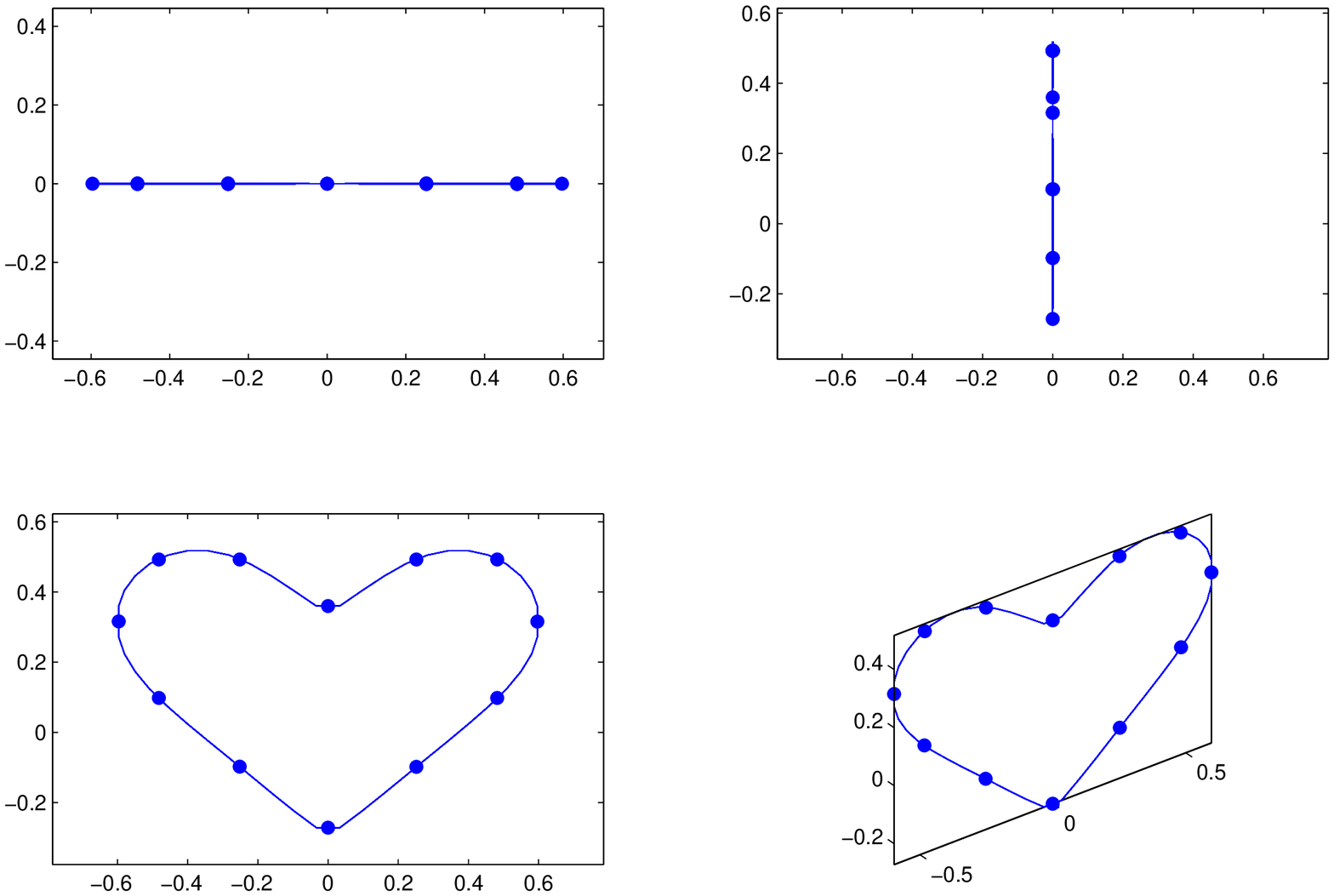}
    \end{minipage}\caption{12 bodies with symmetry 4 for equal masses. (The left figure and right one are in different scales.)}\label{R3X12E01}
\end{figure}

\begin{figure}[!ht]
    \begin{minipage}[t]{0.45\linewidth}
    \includegraphics[width=3in]{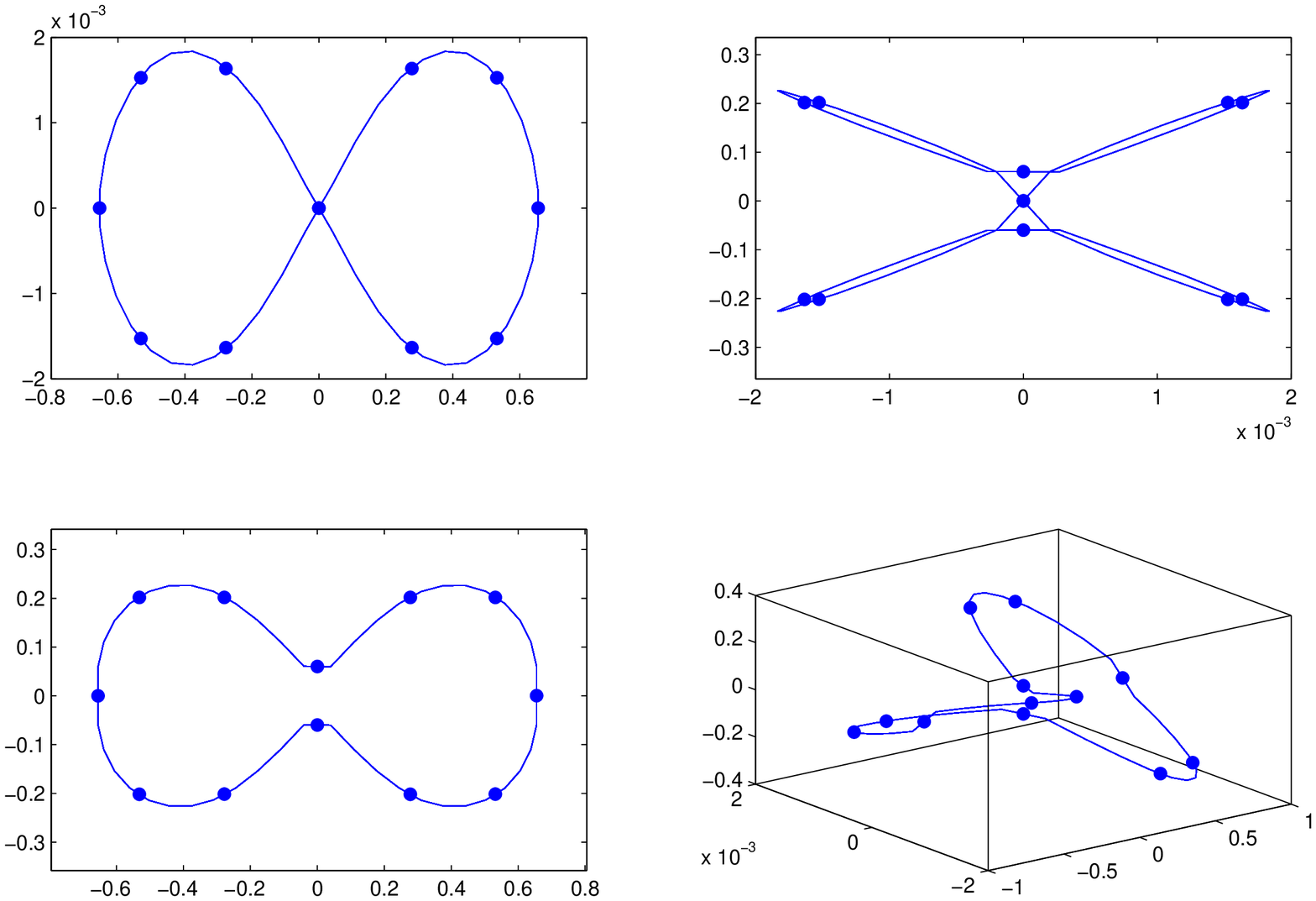}
    \end{minipage}
    \hspace{0.5cm}
    \begin{minipage}[t]{0.45\linewidth}
    \includegraphics[width=3in]{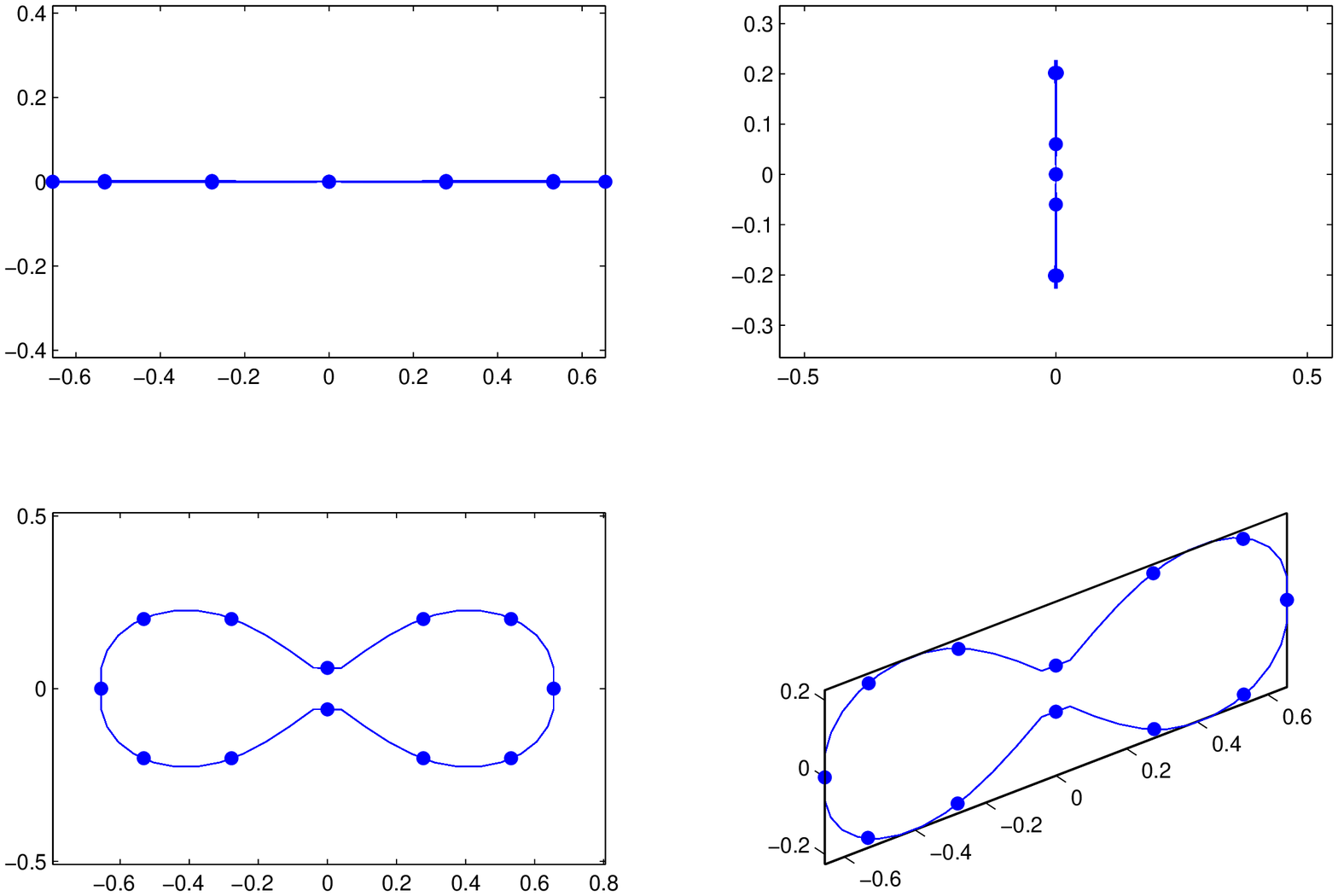}
    \end{minipage}\caption{12 bodies with symmetry 5 for equal masses. (The left figure and right one are in different scales.)}\label{R3X12E02}
\end{figure}

\begin{figure}[!ht]
    \begin{minipage}[t]{0.45\linewidth}
    \includegraphics[width=3in]{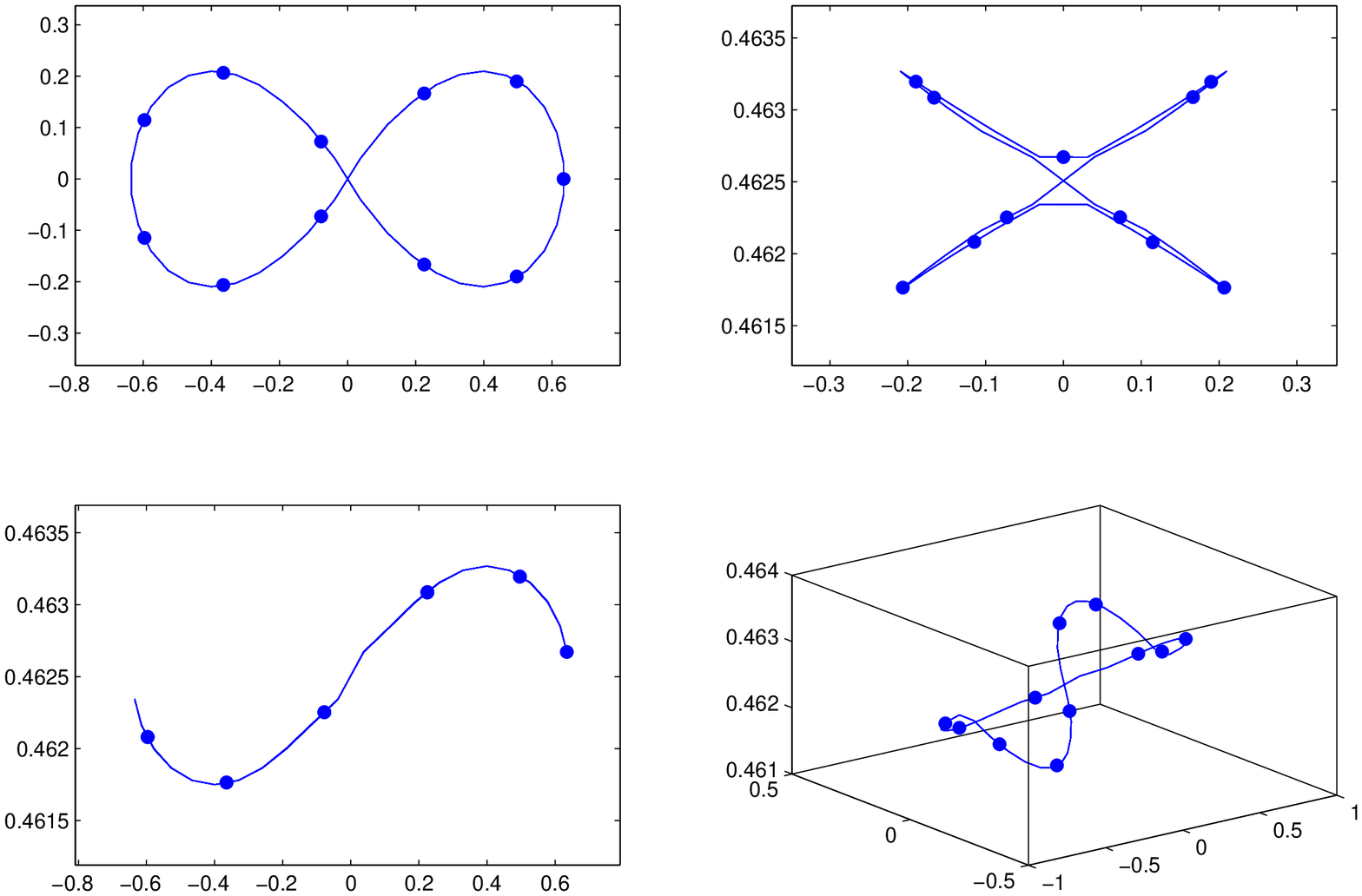}
    \end{minipage}
    \hspace{0.5cm}
    \begin{minipage}[t]{0.45\linewidth}
    \includegraphics[width=3in]{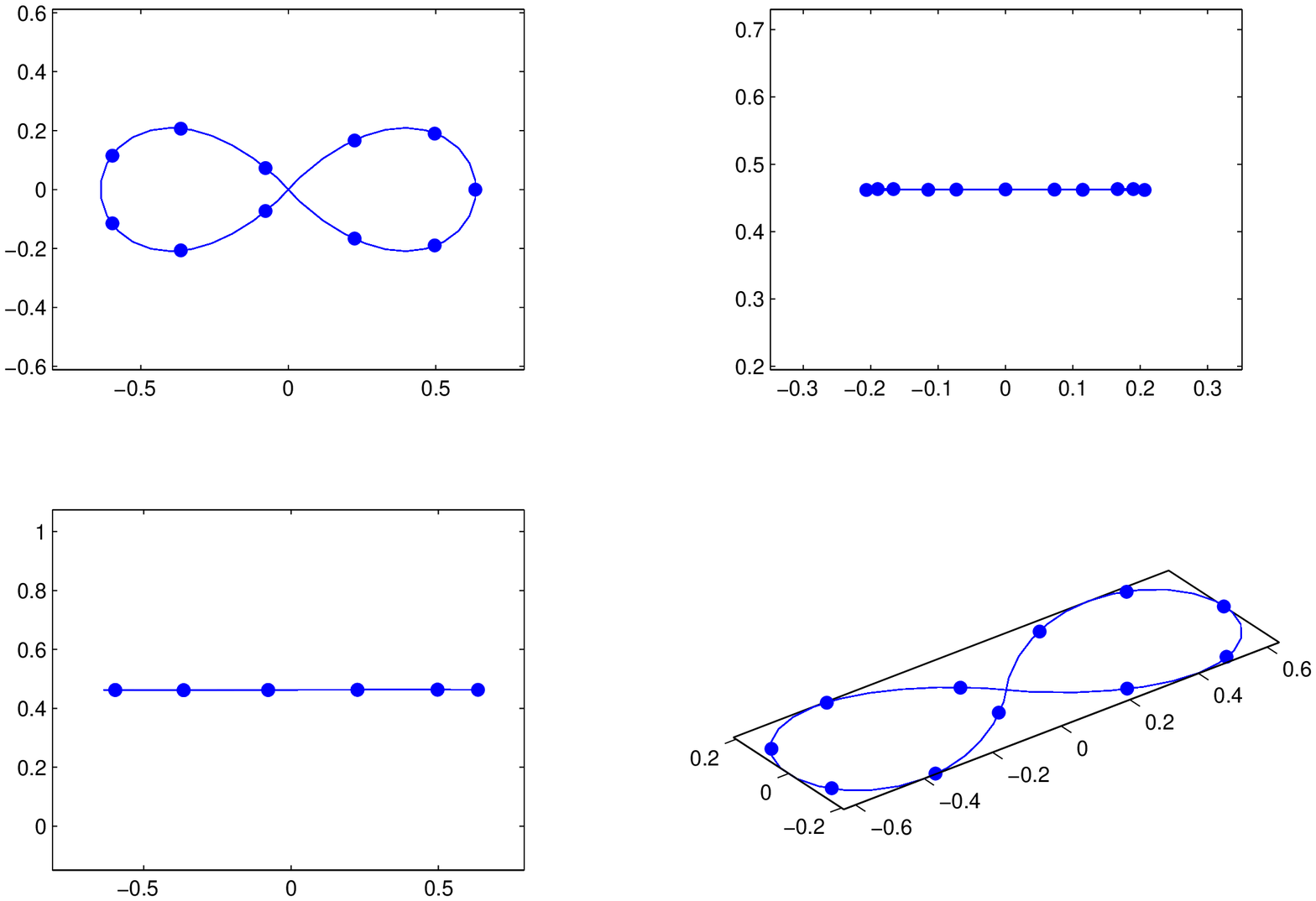}
    \end{minipage}\caption{11 bodies with symmetry 6 for equal masses. (The left figure and right one are in different scales.)}\label{R3X11E03}
\end{figure}

During our numerical study, it shows some interesting phenomena. The
followings are some of them for further consideration:

\begin{itemize}
\item When we set some other constraints on $q^{(3)}$ at will,
e.g. the following symmetries 4-6, we have gotten some interesting
figures numerically (see Figure \ref{R3X12E01}-\ref{R3X11E03}).
\begin{itemize}
\item
symmetry 4
\begin{equation}
q^{(3)}(t)=q^{(3)}\left(t+\frac{T}{4}\right),\quad
t\in\left[0,\frac{T}{4}\right),
\end{equation}
\begin{equation}
q^{(3)}(t)=-q^{(3)}\left(t+\frac{T}{4}\right),\quad
t\in\left[\frac{T}{2},\frac{3T}{4}\right].
\end{equation}

\item
symmetry 5
\begin{equation}
q^{(3)}(t)=q^{(3)}\left(t+\frac{T}{4}\right),\quad
t\in\left[0,\frac{T}{4}\right)\cup\left[\frac{T}{2},\frac{3T}{4}\right].
\end{equation}

\item
symmetry 6
\begin{equation}
q^{(3)}(t)=q^{(3)}\left(t-\frac{3T}{4}\right),\quad
t\in\left[\frac{3T}{4},T\right],
\end{equation}
\begin{equation}
q^{(3)}(t)=q^{(3)}\left(t-\frac{T}{4}\right),\quad
t\in\left[\frac{T}{2},\frac{3T}{4}\right).
\end{equation}
\end{itemize}
See {\bf Figures \ref{R3X12E01}-\ref{R3X11E03}} for numerical
examples.

\item
In Figure-Eight, we can arrange more bodies on it, but we can not
find the ``8'' including more self-intersections, because the dense
subset $I$ is so dense that we can not separate those points from
others by a modern computer.

\item
In Double-Eight, the shape of those two ``8'' depends on the
proportion of the masses. The orbit with less masses looks more
slim, and the bodies on it shuttle rapidly across the interspace of
the bodies on the other orbit.

\item
In Double-Eight, the shape of the orbits also depends on the phase
difference and the initial values of functional $J$.

\item
In symmetry 1, 4, 5 and 6, there are two of the projections of the
orbits become slimmer and slimmer. It seems that one of the
constraints acting on three projections is ``stronger'' than the
other two.

\item
In symmetry 2 and 3, Theorem \ref{th3} guarantees that we can
arrange many bodies on those orbits, no matter even or odd.
\end{itemize}

\begin{figure}[h]
    \begin{minipage}[t]{0.3\linewidth}
    \includegraphics[width=2in]{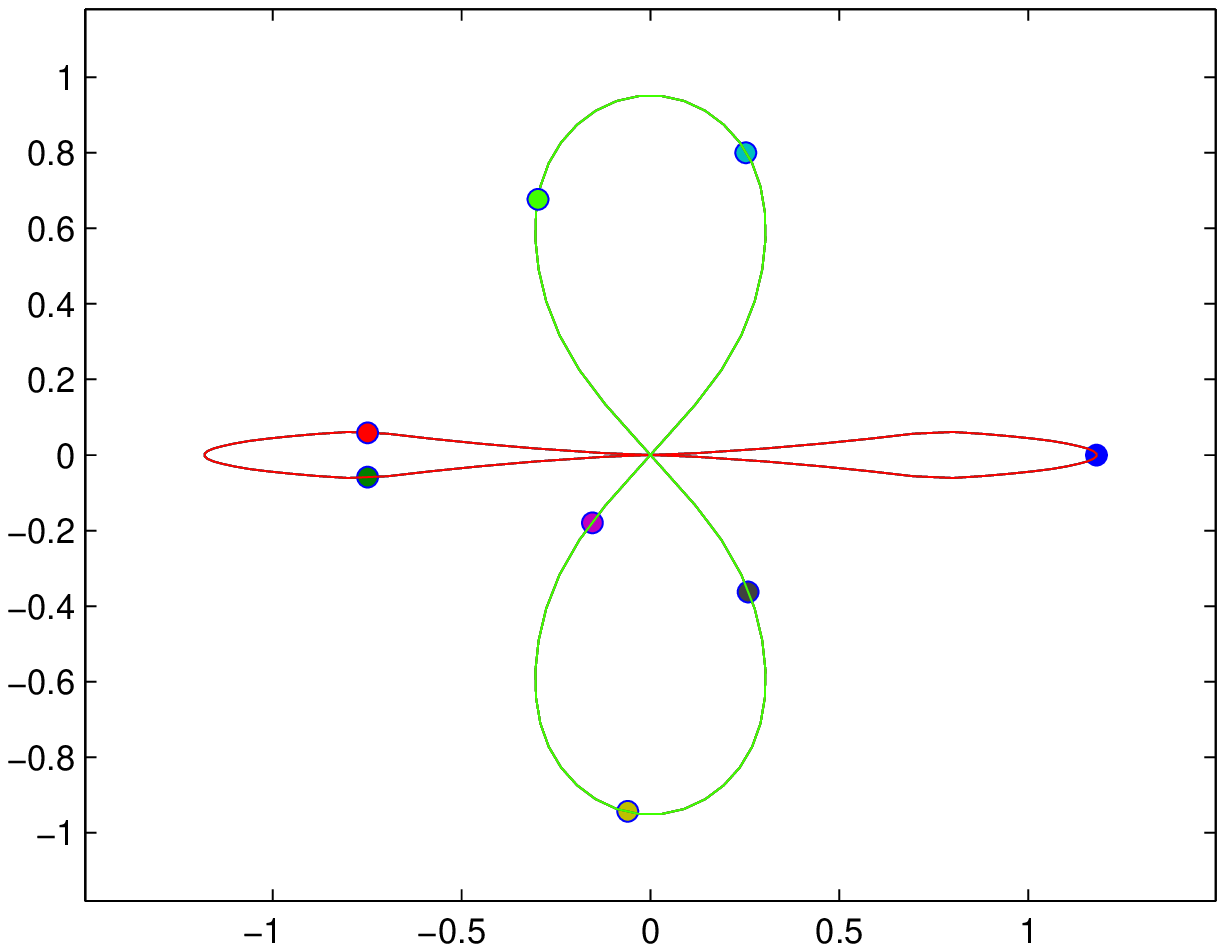}
    \caption{m=[1,1,1;10,10,10,10,10];\newline $\alpha_{1}=0, \alpha_{2}=\frac{T}{6}$.}\label{F3_5N02}
    \end{minipage}
    \hspace{0.3cm}
    \begin{minipage}[t]{0.3\linewidth}
    \includegraphics[width=2in]{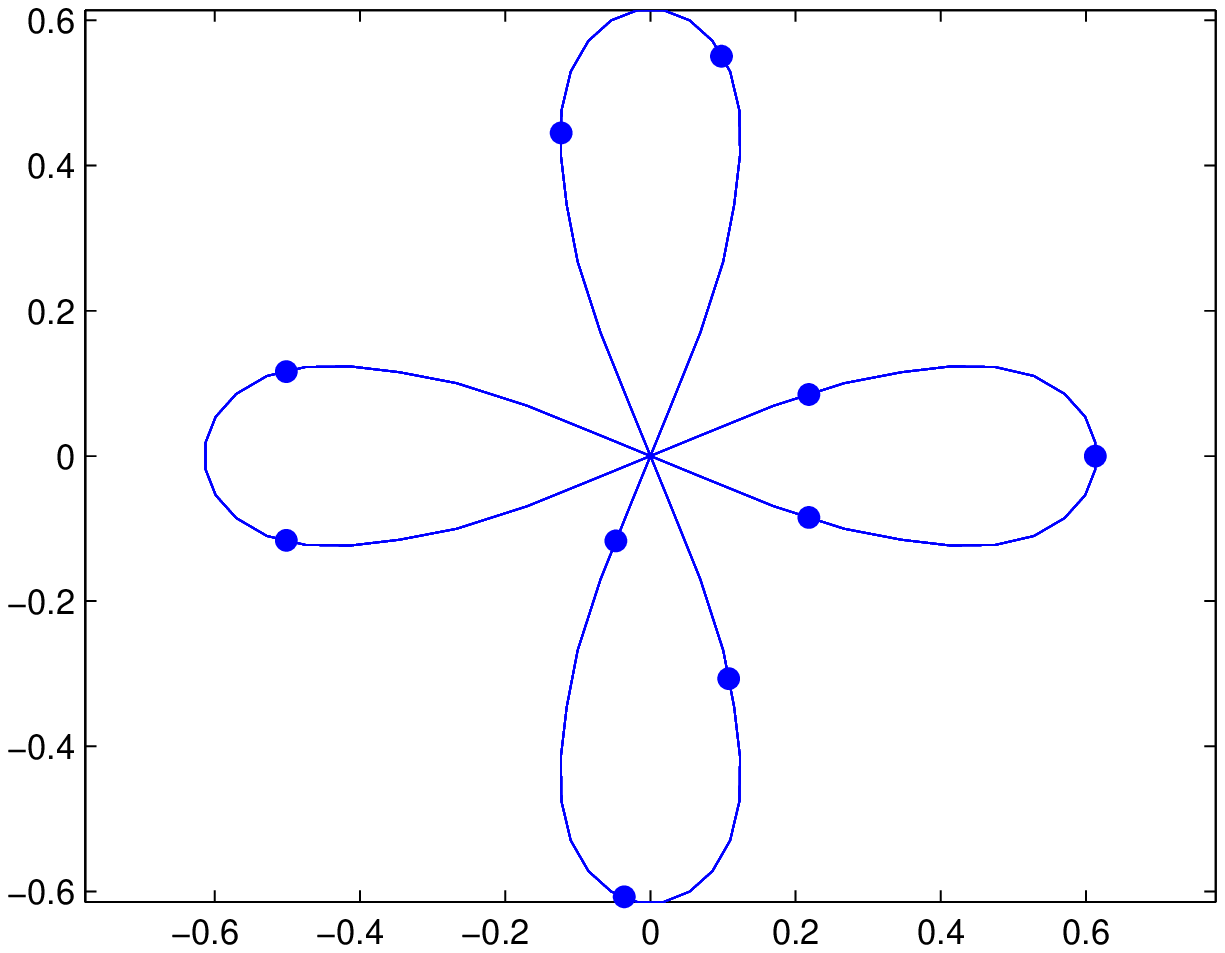}
    \caption{m=[1,1,1,1,1;1,1,1,1,1];\newline $\alpha_{1}=0, \alpha_{2}=\frac{T}{6}$.}\label{F5_5E02}
    \end{minipage}
    \hspace{0.3cm}
    \begin{minipage}[t]{0.3\linewidth}
    \includegraphics[width=2in]{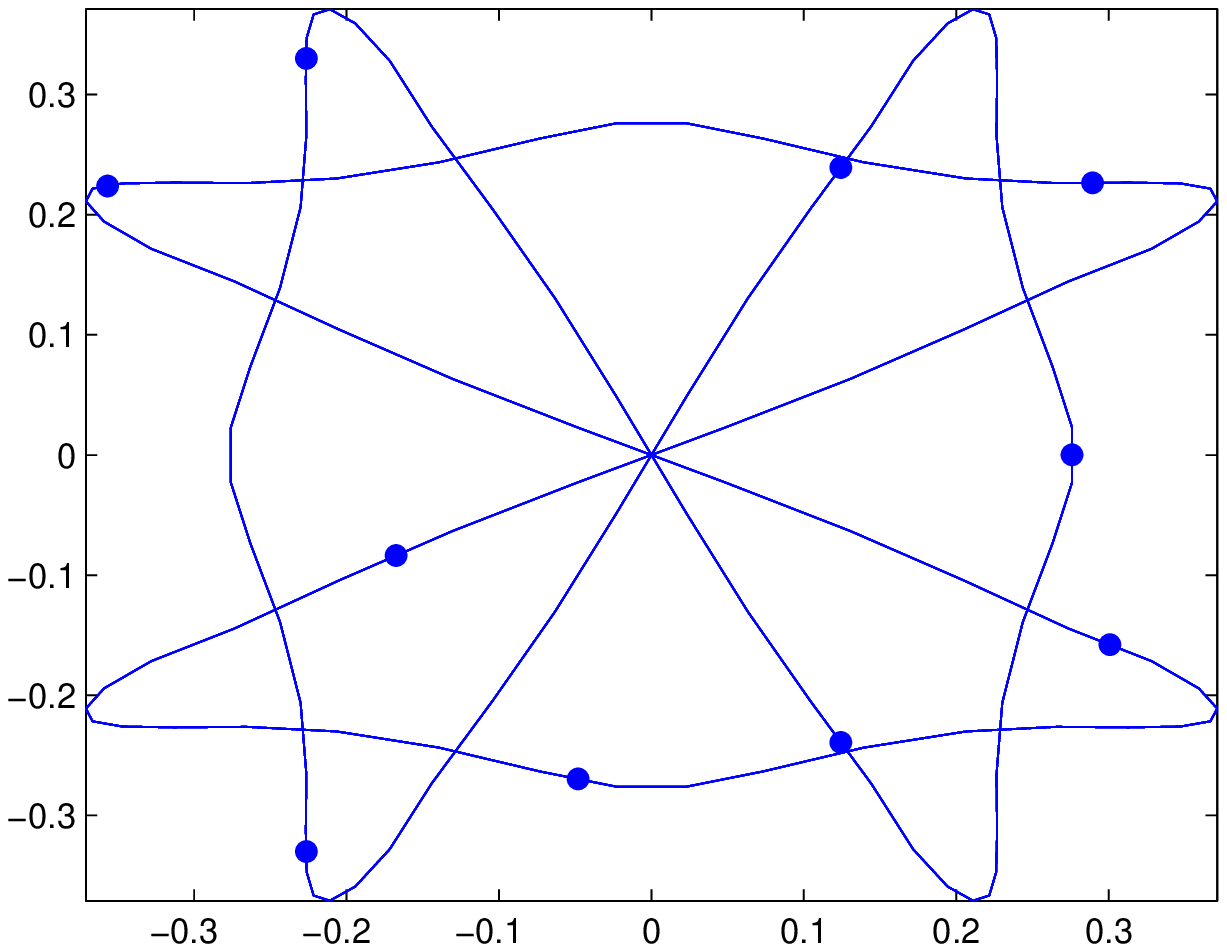}
    \caption{m=[1,1,1,1,1;1,1,1,1,1];\newline $\alpha_{1}=0, \alpha_{2}=\frac{T}{10}$.}\label{F5_5E01}
    \end{minipage}
\end{figure}

\begin{figure}[h]
    \begin{minipage}[t]{0.3\linewidth}
    \includegraphics[width=2in]{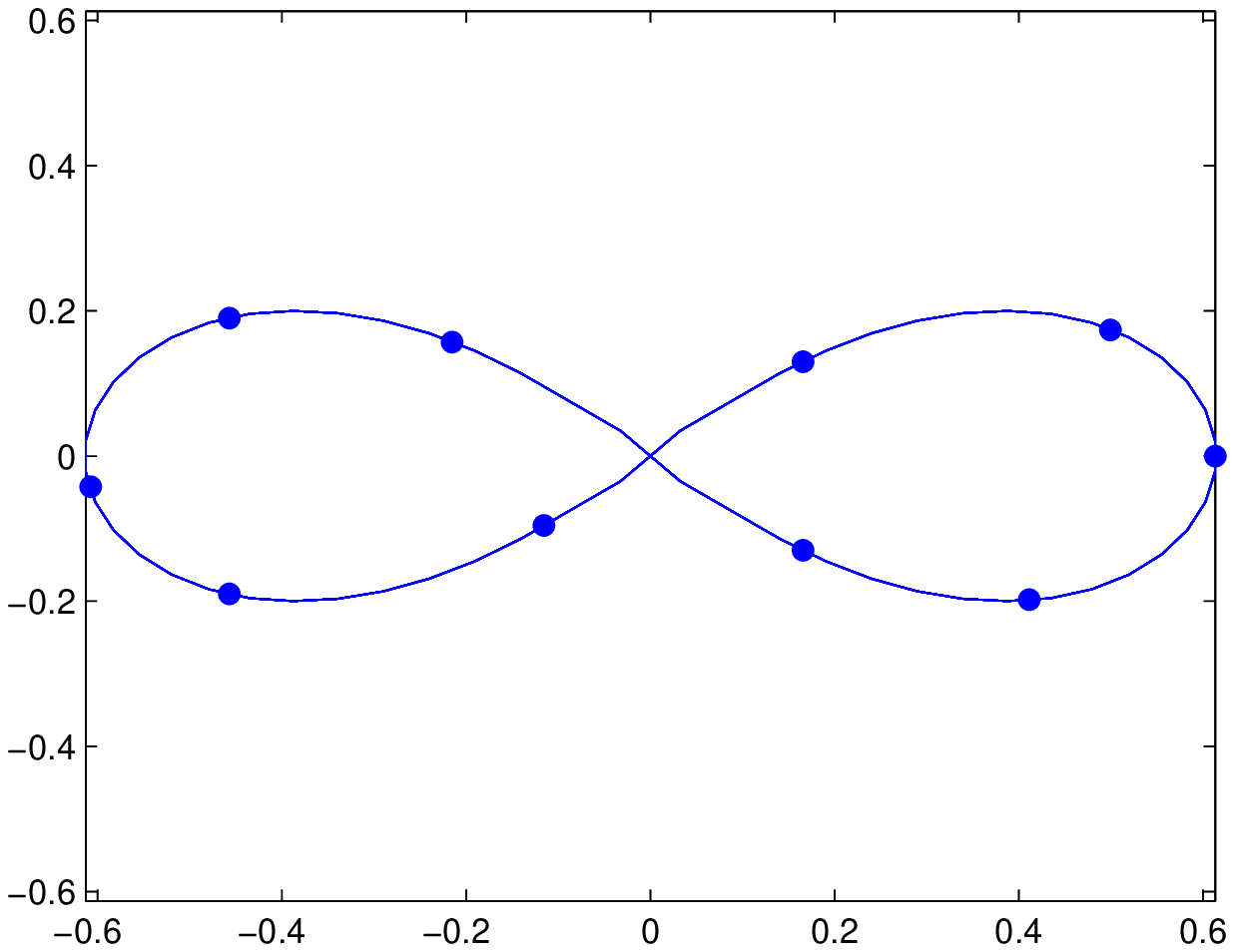}
    \caption{m=[1,1,1,1,1;1,1,1,1,1];\newline $\alpha_{1}=0, \alpha_{2}=\frac{T}{10}$.}\label{L5_5E01}
    \end{minipage}
    \hspace{0.3cm}
    \begin{minipage}[t]{0.3\linewidth}
    \includegraphics[width=2in]{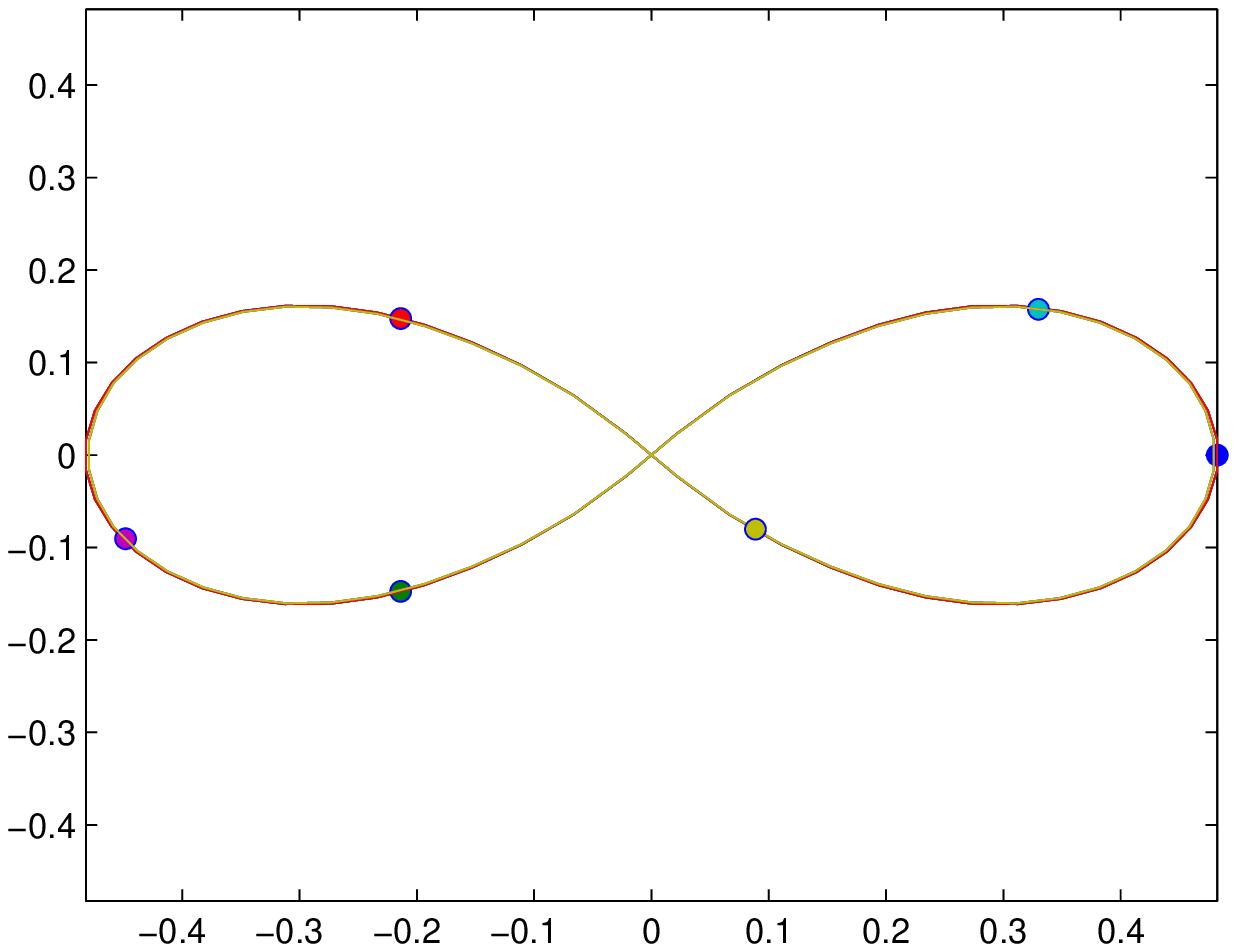}
    \caption{m=[1,1,1;1,1,1];\newline $\alpha_{1}=0, \alpha_{2}=\frac{T}{8}$.}\label{L3_3E04}
    \end{minipage}
    \hspace{0.3cm}
    \begin{minipage}[t]{0.3\linewidth}
    \includegraphics[width=2in]{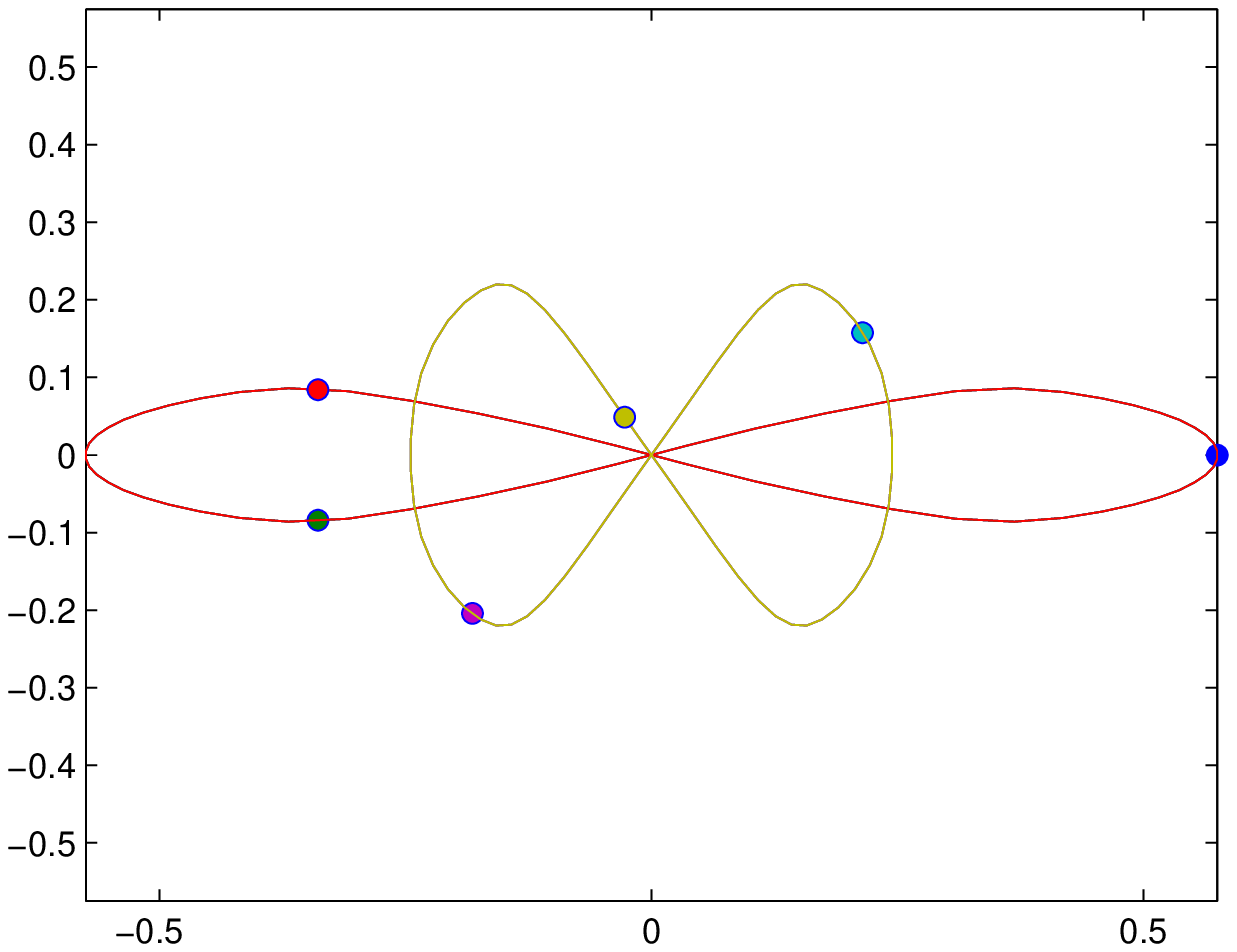}
    \caption{m=[1,1,1;1,1,1];\newline $\alpha_{1}=0, \alpha_{2}=\frac{T}{12}$.}\label{L3_3E05}
    \end{minipage}
\end{figure}

\begin{figure}[h]
    \begin{minipage}[t]{0.3\linewidth}
    \includegraphics[width=2in]{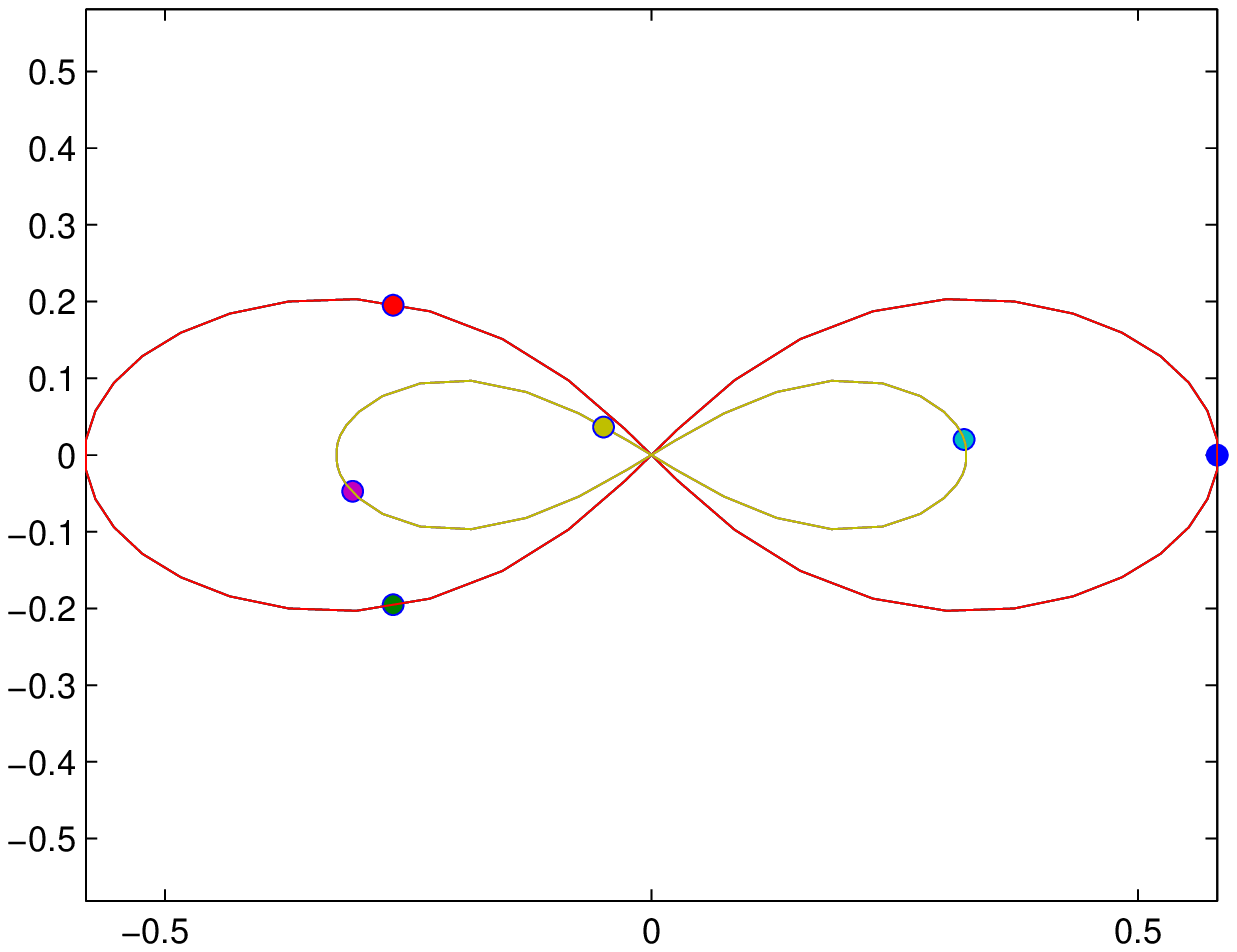}
    \caption{m=[1,1,1;1,1,1];\newline $\alpha_{1}=0, \alpha_{2}=\frac{T}{12}$.}\label{L3_3E03}
    \end{minipage}
    \hspace{0.3cm}
    \begin{minipage}[t]{0.3\linewidth}
    \includegraphics[width=2in]{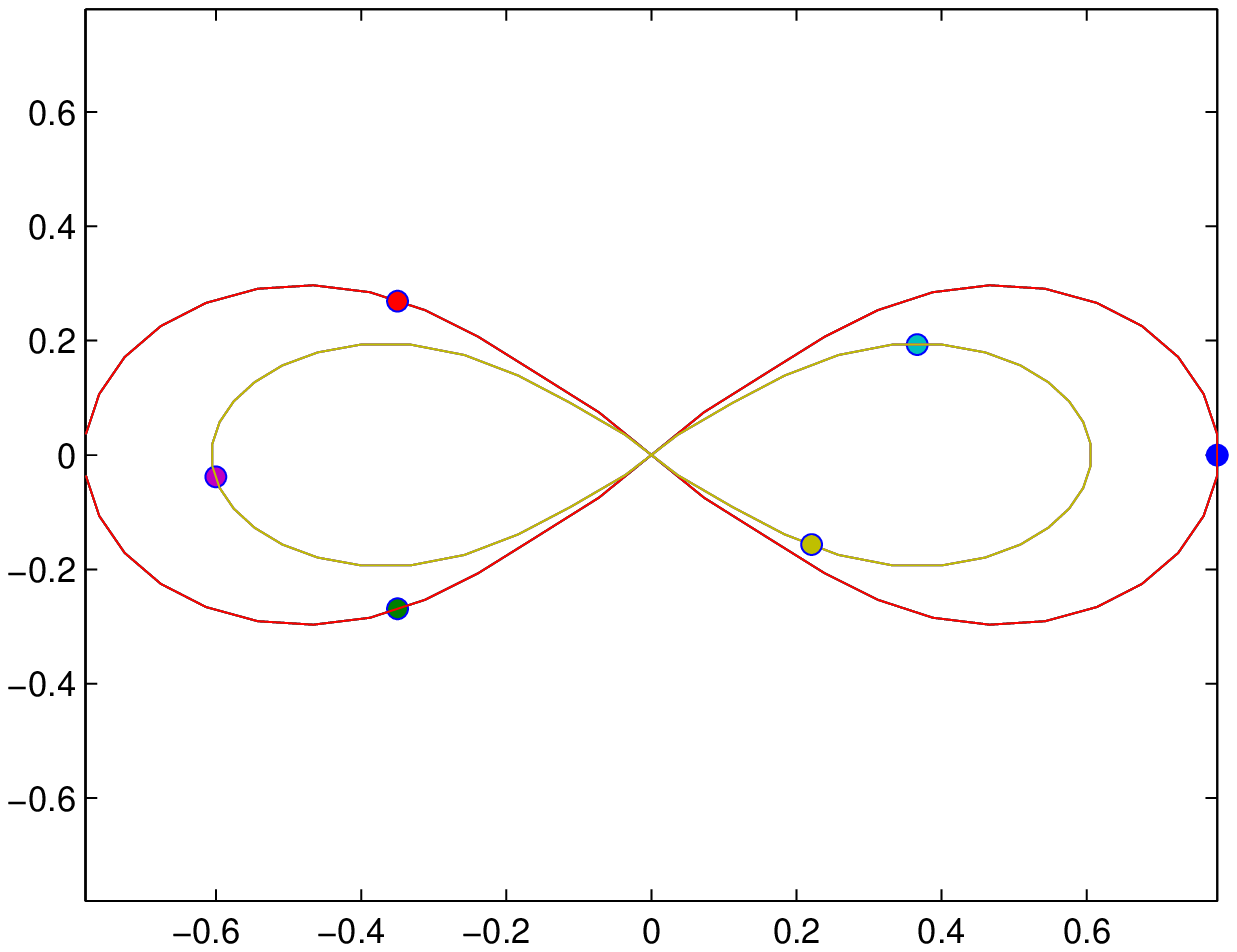}
    \caption{m=[1,1,1;5,5,5];\newline $\alpha_{1}=0, \alpha_{2}=\frac{T}{6}$.}\label{L3_3N01}
    \end{minipage}
    \hspace{0.3cm}
    \begin{minipage}[t]{0.3\linewidth}
    \includegraphics[width=2in]{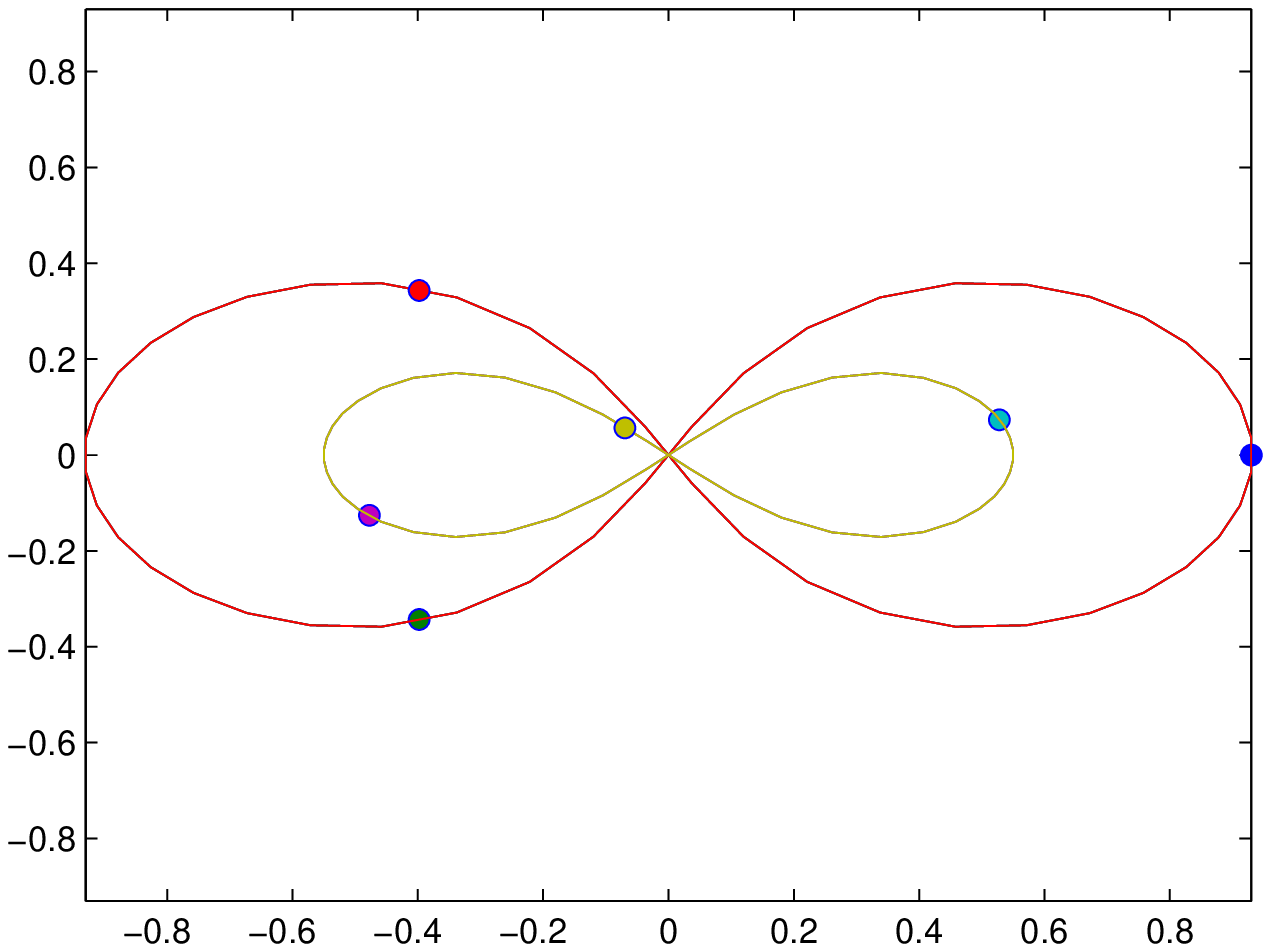}
    \caption{m=[1,1,1;5,5,5];\newline $\alpha_{1}=0, \alpha_{2}=\frac{T}{12}$.}\label{L3_3N02}
    \end{minipage}
\end{figure}

\begin{figure}[h]
    \begin{minipage}[t]{0.3\linewidth}
    \includegraphics[width=2in]{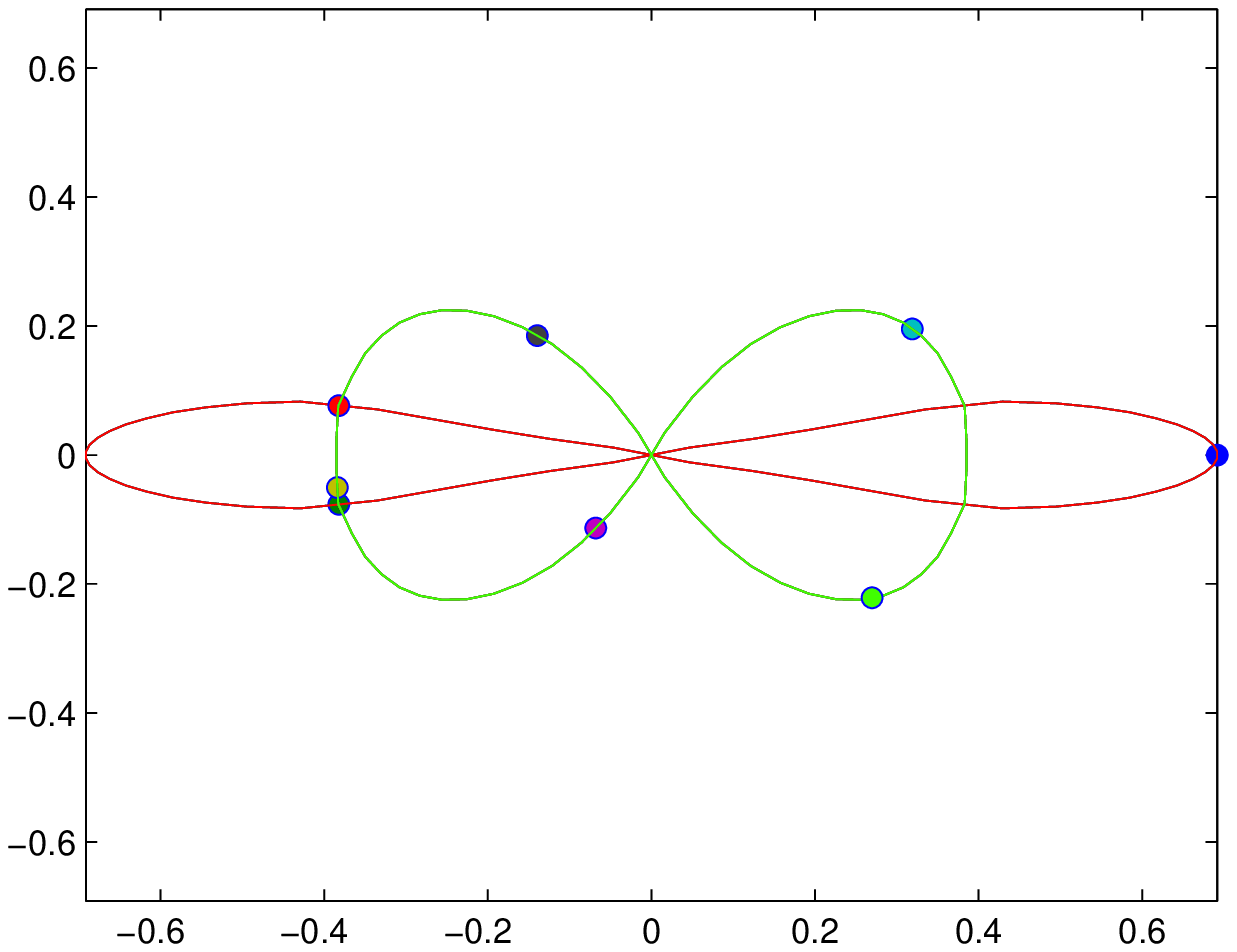}
    \caption{m=[1,1,1;1,1,1,1,1];\newline $\alpha_{1}=0, \alpha_{2}=\frac{T}{10}$.}\label{L3_5E01}
    \end{minipage}
    \hspace{0.3cm}
    \begin{minipage}[t]{0.3\linewidth}
    \includegraphics[width=2in]{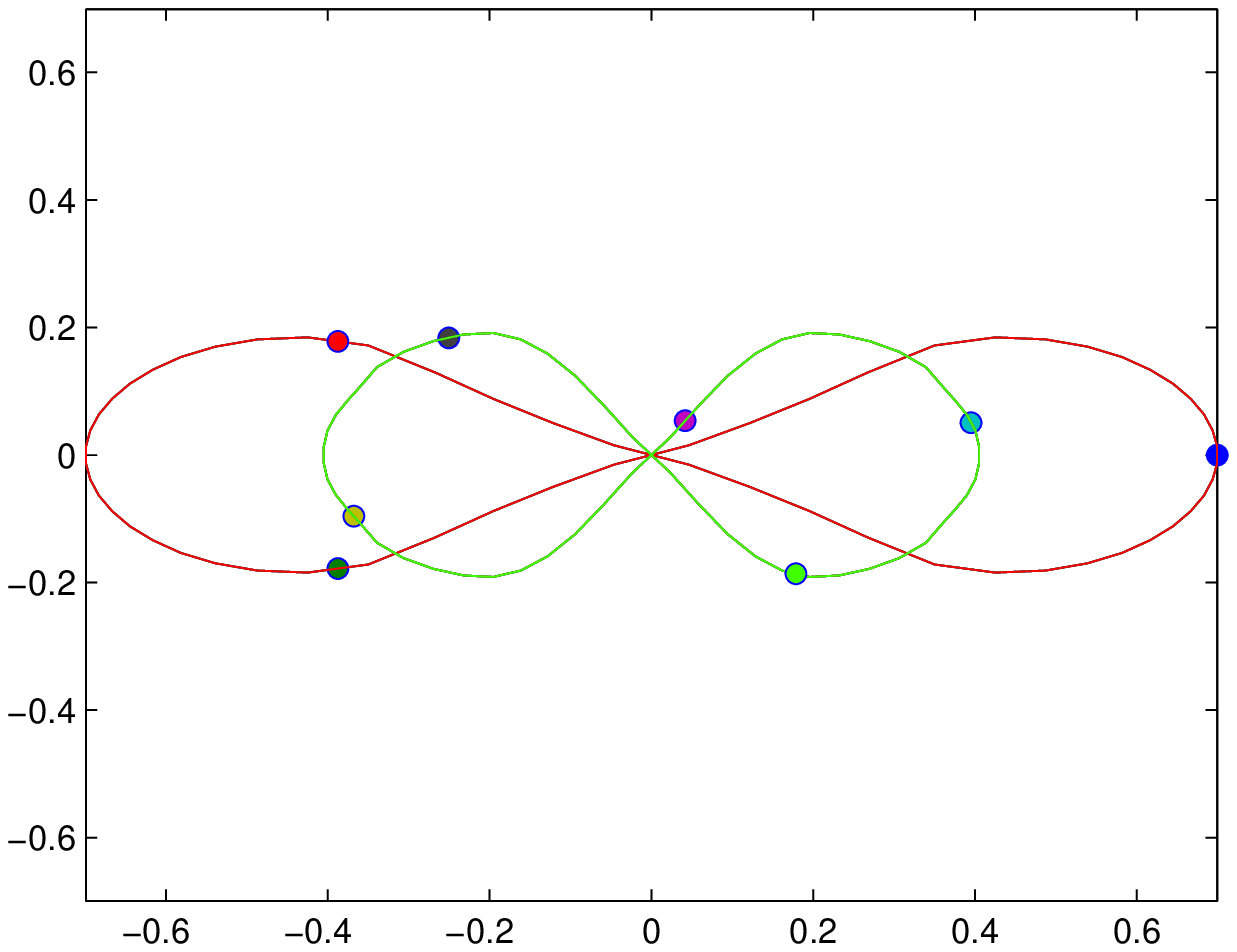}
    \caption{m=[1,1,1;1,1,1,1,1];\newline $\alpha_{1}=0, \alpha_{2}=\frac{T}{20}$.}\label{L3_5E02}
    \end{minipage}
    \hspace{0.3cm}
    \begin{minipage}[t]{0.3\linewidth}
    \includegraphics[width=2in]{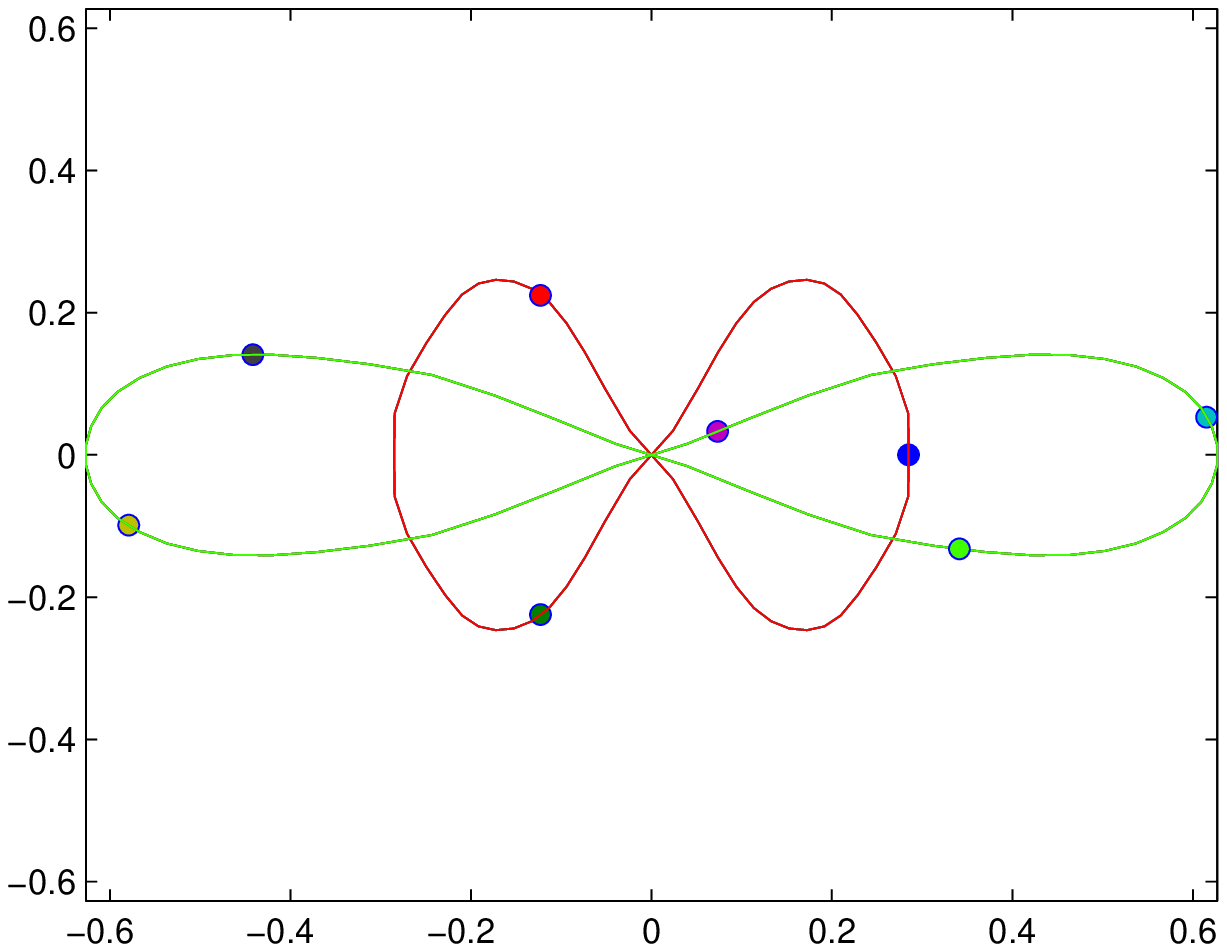}
    \caption{m=[1,1,1;1,1,1,1,1];\newline $\alpha_{1}=0, \alpha_{2}=\frac{T}{20}$.}\label{L3_5E03}
    \end{minipage}
\end{figure}

\begin{figure}[h]
    \begin{minipage}[t]{0.3\linewidth}
    \includegraphics[width=2in]{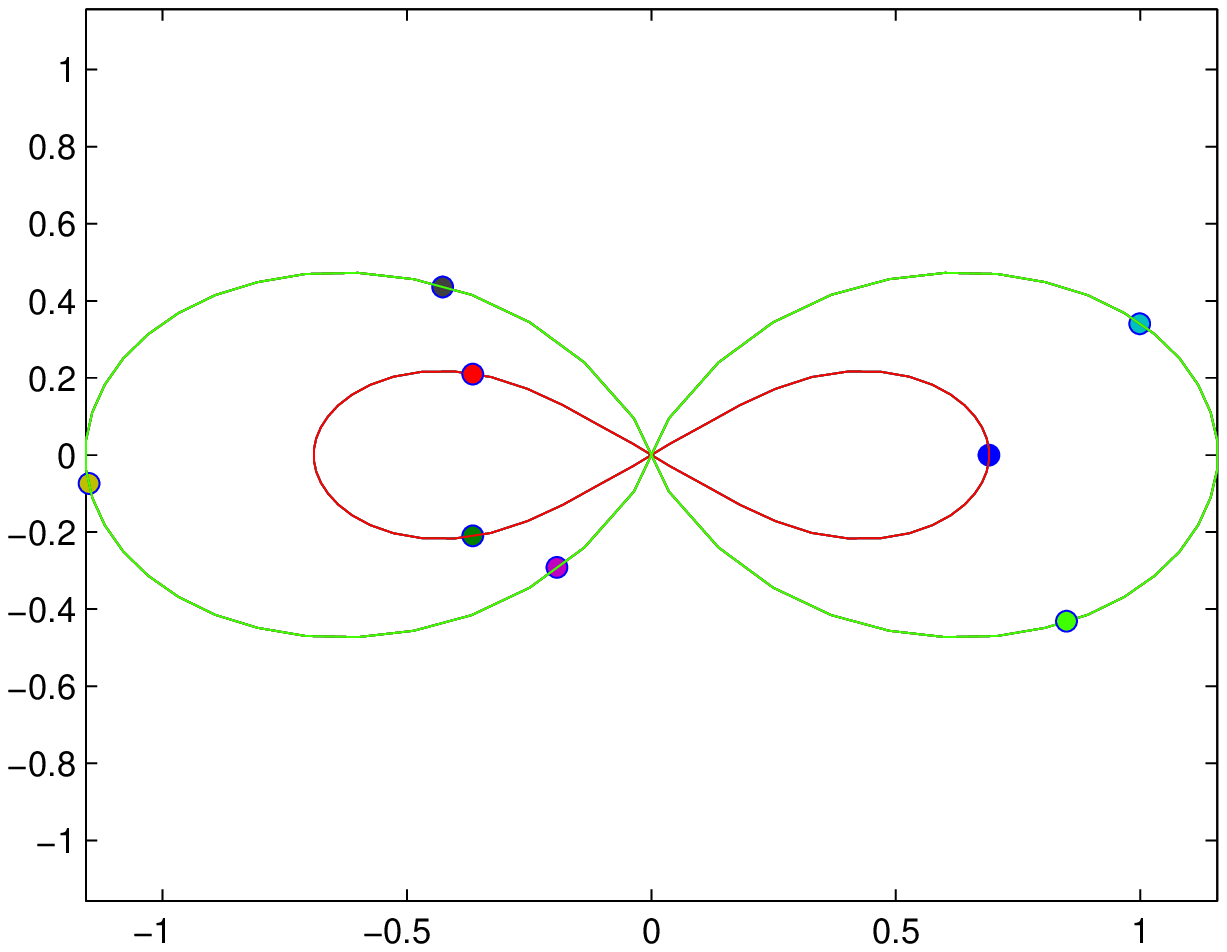}
    \caption{m=[10,10,10;1,1,1,1,1];\newline $\alpha_{1}=0, \alpha_{2}=\frac{T}{10}$.}\label{L3_5N09}
    \end{minipage}
    \hspace{0.3cm}
    \begin{minipage}[t]{0.3\linewidth}
    \includegraphics[width=2in]{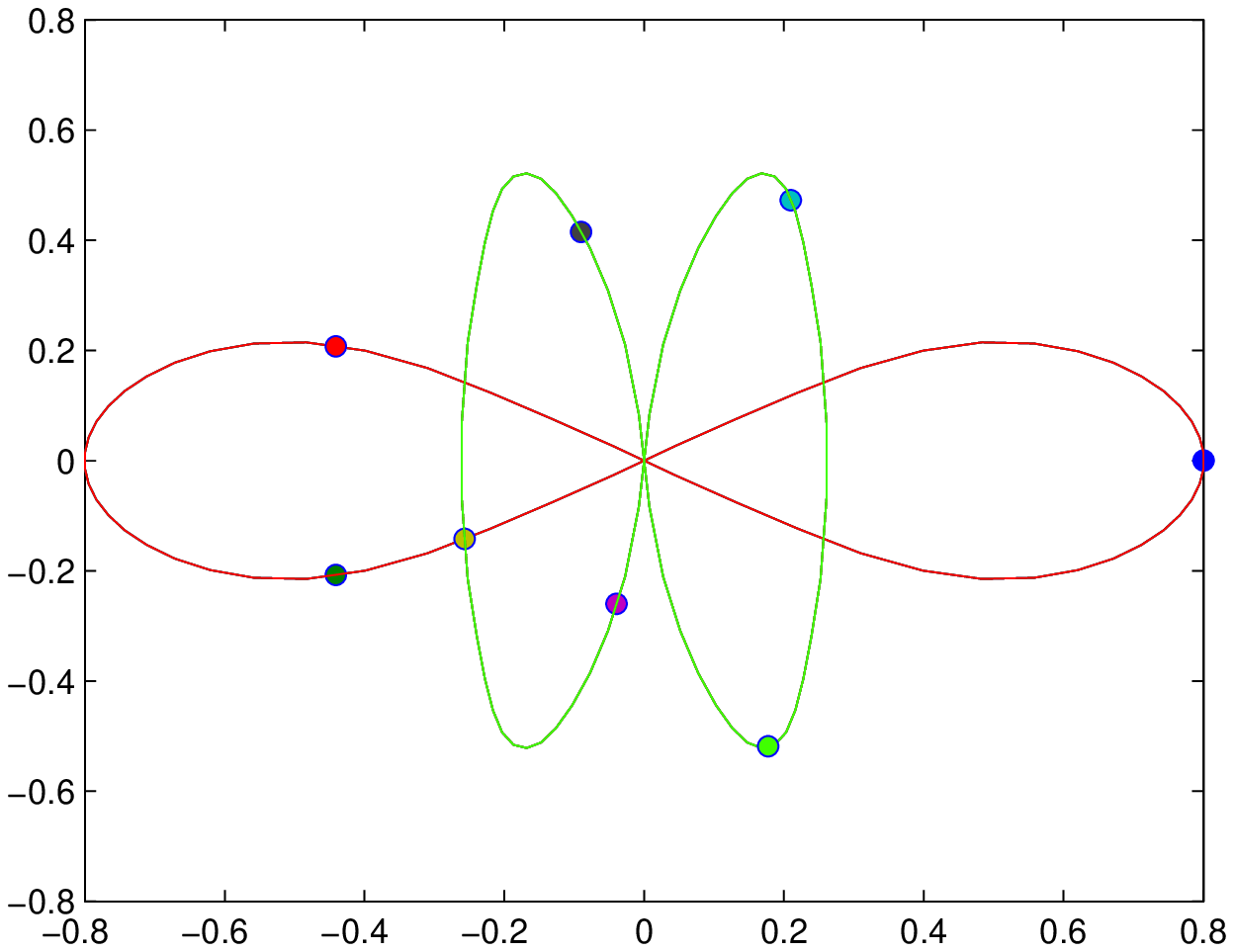}
    \caption{m=[10,10,10;1,1,1,1,1];\newline $\alpha_{1}=0, \alpha_{2}=\frac{T}{10}$.}\label{L3_5N03}
    \end{minipage}
    \hspace{0.3cm}
    \begin{minipage}[t]{0.3\linewidth}
    \includegraphics[width=2in]{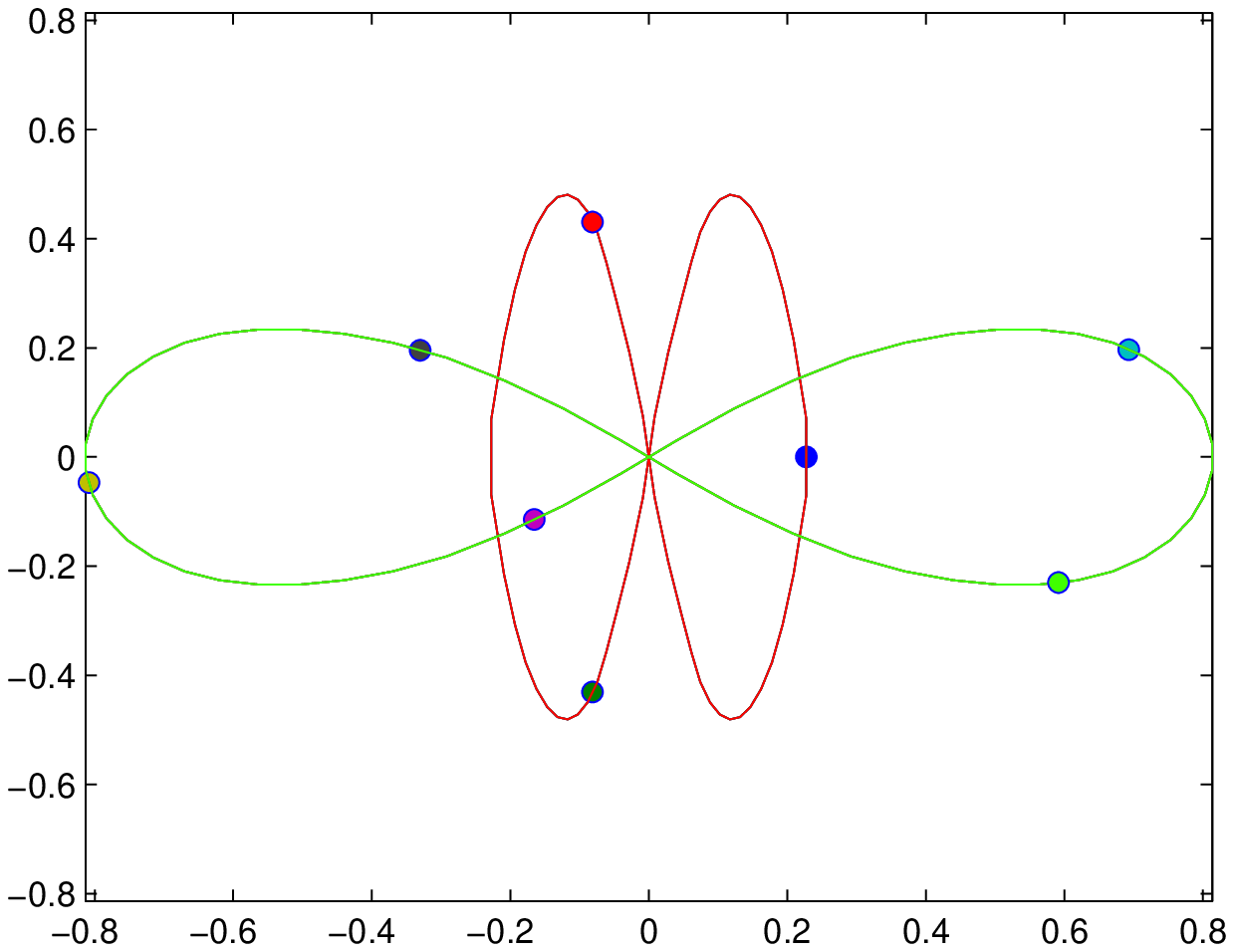}
    \caption{m=[1,1,1;5,5,5,5,5];\newline $\alpha_{1}=0, \alpha_{2}=\frac{T}{10}$.}\label{L3_5N08}
    \end{minipage}
\end{figure}

The following figures are numerical examples for Theorem \ref{th3}:

\begin{figure}[h]
    \begin{minipage}[t]{0.3\linewidth}
    \includegraphics[width=2in]{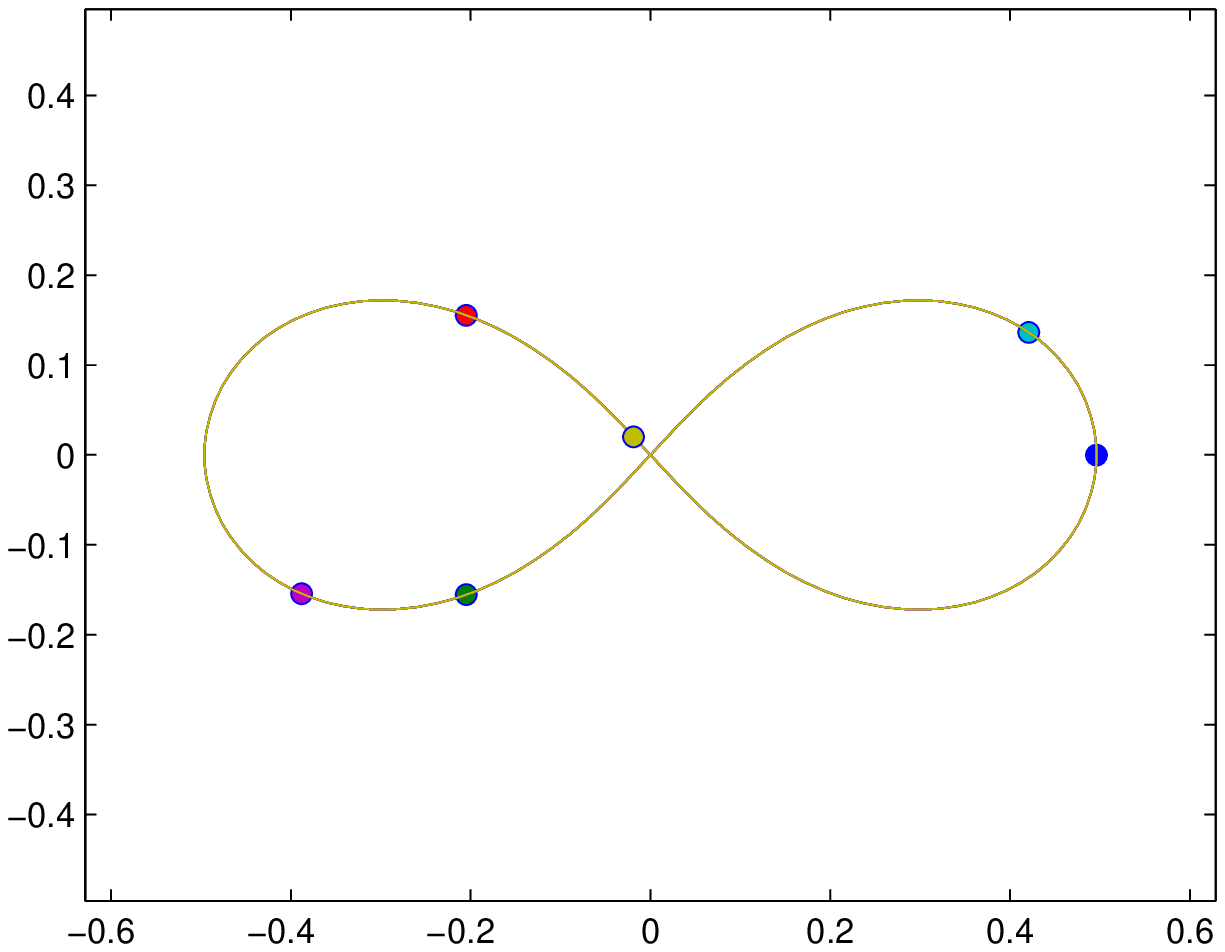}
    \caption{6 bodies on Figure-Eight:\newline $m=[1,1,1;1,1,1],\newline \alpha_{1}=0, \alpha_{2}=\frac{T}{12}$.}\label{Even06E01}
    \end{minipage}
    \hspace{0.3cm}
    \begin{minipage}[t]{0.3\linewidth}
    \includegraphics[width=2in]{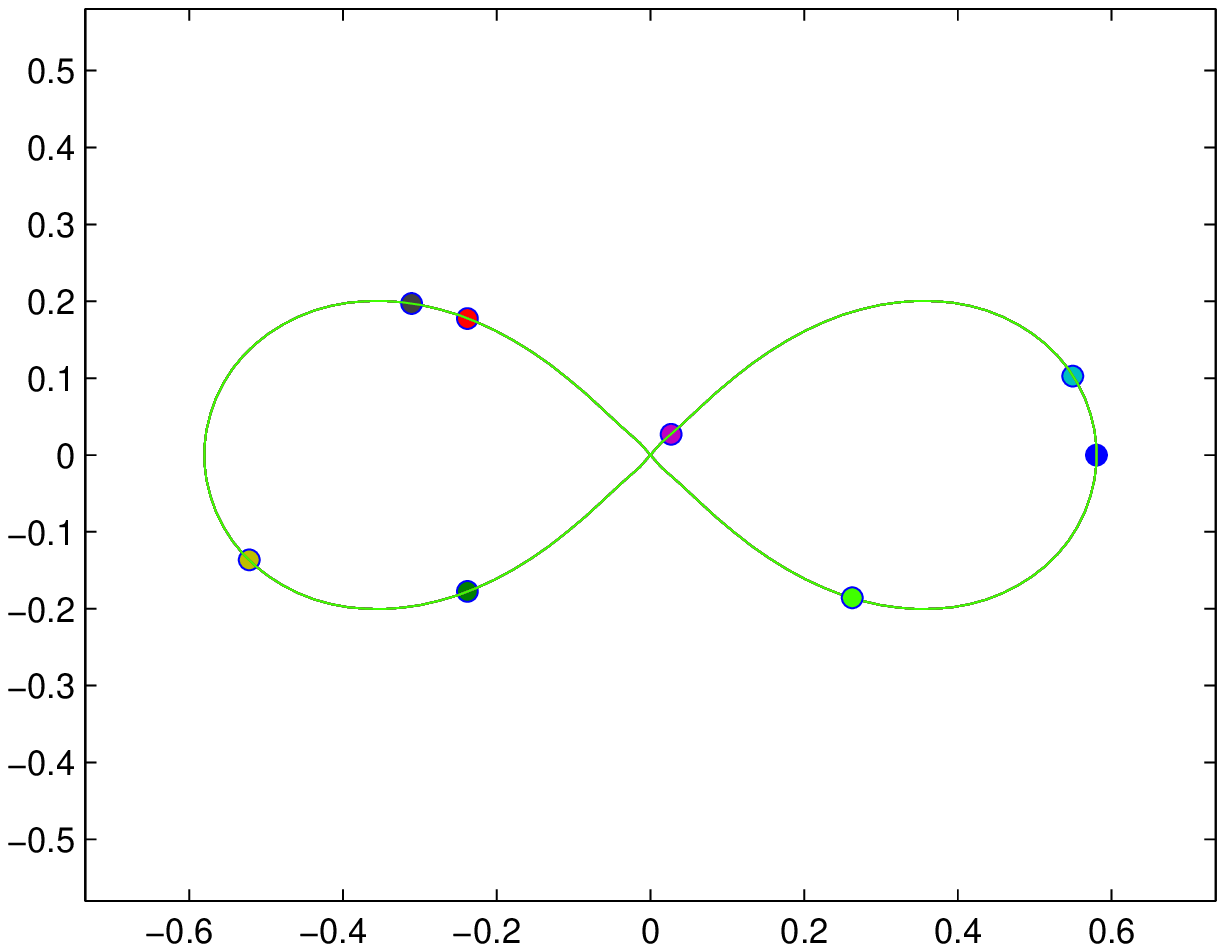}
    \caption{8 bodies on Figure-Eight: $m=[1,1,1;1,1,1,1,1],\newline \alpha_{1}=0, \alpha_{2}=\frac{T}{20}$.}\label{Even08E01}
    \end{minipage}
    \hspace{0.3cm}
    \begin{minipage}[t]{0.3\linewidth}
    \includegraphics[width=2in]{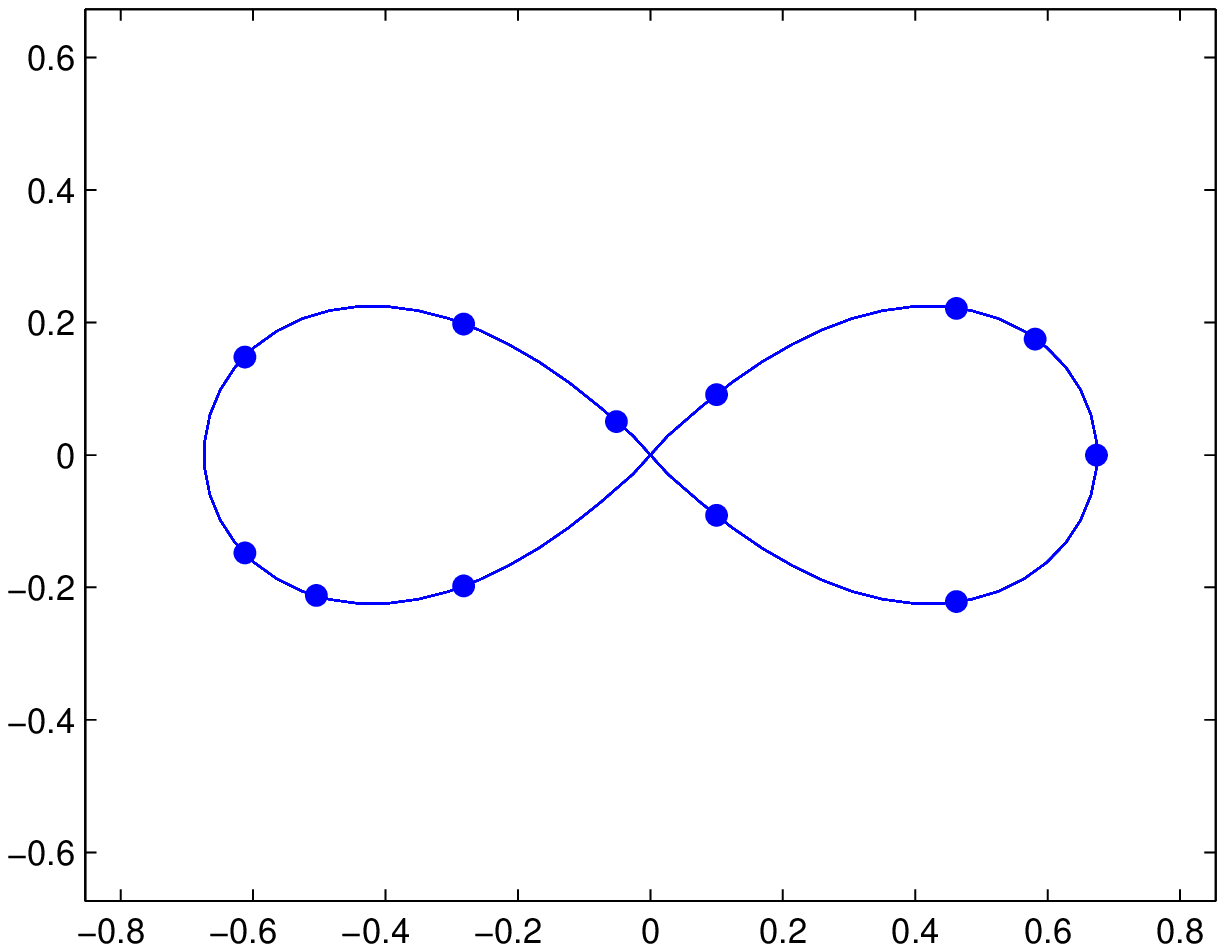}
    \caption{12 bodies on Figure-Eight: $m=[1,1,1;1,1,1,1,1,1,1,1,1],\newline  \alpha_{1}=0, \alpha_{2}=\frac{T}{36}$.}\label{Even12E01}
    \end{minipage}
\end{figure}

\clearpage
\appendix
\section{The estimation of $|q|$ for other symmetric constraints}
{\em For Double-Eight}:
\begin{eqnarray*}
\left|q(t)\right| & \leq &
\left|q^{(1)}(t)\right|+\left|q^{(2)}(t)\right|\\
& =  & \frac{1}{2}\left|q^{(1)}(t)-q^{(1)}(t+\frac{T}{2})\right|
+\frac{1}{2}\left|q^{(2)}(t)-q^{(2)}(t+\frac{T}{4})\right|\\
& =  &
\frac{1}{2}\left|\int_{t}^{t+\frac{T}{2}}\dot{q}^{(1)}(s)ds\right|
+\frac{1}{2}\left|\int_{t}^{t+\frac{T}{4}}\dot{q}^{(2)}(s)ds\right|\\
&\leq&
\frac{1}{2}\int_{t}^{t+\frac{T}{2}}\left|\dot{q}^{(1)}(s)\right|ds
+\frac{1}{2}\int_{t}^{t+\frac{T}{4}}\left|\dot{q}^{(2)}(s)\right|ds\\
&\leq&
\frac{1}{2}\left(\int_{t}^{t+\frac{T}{2}}\left|\dot{q}^{(1)}(s)\right|^{2}ds\right)^{\frac{1}{2}}
\left(\frac{T}{2}\right)^{\frac{1}{2}}\\
&    &
+\frac{1}{2}\left(\int_{t}^{t+\frac{T}{4}}\left|\dot{q}^{(2)}(s)\right|^{2}ds\right)^{\frac{1}{2}}
\left(\frac{T}{4}\right)^{\frac{1}{2}}\\
&\leq&
\frac{\sqrt{T}}{4}\left(\left(\int_{0}^{T}\left|\dot{q}^{(1)}(t)\right|^{2}dt\right)^{\frac{1}{2}}
+\left(\int_{0}^{T}\left|\dot{q}^{(2)}(t)\right|^{2}dt\right)^{\frac{1}{2}}\right)\\
&\leq&
\frac{\sqrt{3T}}{4}\left(\int_{0}^{T}\left|\dot{q}^{(1)}(t)\right|^{2}dt
+\int_{0}^{T}\left|\dot{q}^{(2)}(t)\right|^{2}dt\right)^{\frac{1}{2}}\\
&=&
\frac{\sqrt{3T}}{4}\left(\int_{0}^{T}\left|\dot{q}(t)\right|^{2}dt\right)^{\frac{1}{2}}.
\end{eqnarray*}
For $Q(t)$ in Double-Eight, we have the same estimation.
\newline
{\em For symmetry 2}: Since (\ref{con2}) is the same condition as in
symmetry 1, we can at least draw the same conclusion as in symmetry
1, i.e.
\begin{equation*}
\left|q(t)\right| \leq
\frac{\sqrt{3T}}{4}\left(\int_{0}^{T}\left|\dot{q}(t)\right|^{2}dt\right)^{\frac{1}{2}}.
\end{equation*}\ \
\newline
{\em For symmetry 3}: We set $X=(q,Q)$, it is easy to verify that
\begin{equation*}
\left|X(t)\right| \leq\left|q(t)\right|+\left|Q(t)\right|\leq
\frac{\sqrt{3T}}{2}\left(\int_{0}^{T}\left|\dot{q}(t)\right|^{2}dt\right)^{\frac{1}{2}}.
\end{equation*}\ \

\end{document}